\documentclass[11pt]{article}

\usepackage[T1]{fontenc}
\usepackage[utf8]{inputenc}
\usepackage{lmodern}
\usepackage[a4paper,margin=1in]{geometry}
\usepackage{amsmath,amssymb,amsthm,mathtools}
\usepackage{enumitem}
\usepackage{chapterbib}
\usepackage[hidelinks]{hyperref}
\usepackage{xcolor}

\numberwithin{equation}{section}
\newtheorem{theorem}{Theorem}[section]
\newtheorem{proposition}[theorem]{Proposition}
\newtheorem{lemma}[theorem]{Lemma}
\newtheorem{corollary}[theorem]{Corollary}
\theoremstyle{definition}
\newtheorem{definition}[theorem]{Definition}
\theoremstyle{remark}

\newcommand{\F}{\mathbb F}
\newcommand{\Q}{\mathbb Q}
\newcommand{\Z}{\mathbb Z}
\newcommand{\N}{\mathbb N}
\newcommand{\cO}{\mathcal O}

\newcommand{\fm}{\mathfrak m}

\DeclareMathOperator{\GL}{GL}
\DeclareMathOperator{\SL}{SL}
\DeclareMathOperator{\Nrd}{Nrd}
\DeclareMathOperator{\Trd}{Trd}
\DeclareMathOperator{\gr}{gr}
\DeclareMathOperator{\Ann}{Ann}
\DeclareMathOperator{\Hom}{Hom}
\DeclareMathOperator{\Ext}{Ext}
\DeclareMathOperator{\Supp}{Supp}
\DeclareMathOperator{\Spec}{Spec}
\DeclareMathOperator{\Frac}{Frac}

\DeclareMathOperator{\GKdim}{GKdim}
\DeclareMathOperator{\Krdim}{Kr.dim}
\DeclareMathOperator{\Der}{Der}

\title{Ramified tori and Gelfand--Kirillov dimension\\
of mod $p$ representations}
\author{Jiusi Zhao\thanks{Academy of Mathematics and Systems Science,
Chinese Academy of Sciences, Beijing 100190, China.
\textit{E-mail:}
\href{mailto:zhaojiusi@amss.ac.cn}{\texttt{zhaojiusi@amss.ac.cn}}.}\\[0.4em]
\small with an appendix by Yongquan Hu}
\date{}

\begin{document}

\maketitle

\begin{abstract}
Let $K/\Q_p$ be a finite extension of degree $n$ with $p\geq3$,
and let $D$ be the non-split quaternion algebra over $K$. 
We prove that every infinite-dimensional admissible smooth
representation of $\GL_2(K)$ or $D^\times$ with a central character
and satisfying a Casimir annihilation condition has
Gelfand--Kirillov dimension $n$. As a corollary, every infinite-dimensional admissible Banach
representation whose locally analytic vectors have an infinitesimal
character has dimension $n$, strengthening Sorgdrager's bound to an
equality. An appendix by Yongquan Hu proves that, when $K\ne\Q_p$
and the coefficient field is algebraically closed, the canonical
diagram of an admissible supersingular representation satisfying
the same Casimir annihilation condition is $(\pi,\pi,\mathrm{id})$, 
as in the case of local fields of characteristic $p$.
\end{abstract}

\section{Introduction}\label{sec:introduction}

Let $p\geq3$ be a prime number and let $K/\Q_p$ be a finite
extension, with ramification index $e$, residue degree $f$, and
$n=[K:\Q_p]=ef$. Let $D$ be the non-split quaternion algebra
over $K$, and let $\F$ be a field of characteristic $p$ 
with an embedding $\mathbf F_{p^{2f}}\hookrightarrow\F$.
The mod-$p$ Langlands program studies smooth
admissible representations of $\GL_2(K)$ arising from the
cohomology of Shimura varieties. The predicted mod-$p$
Jacquet--Langlands correspondence also leads to the study of
smooth admissible representations of $D^\times$. 
The Gelfand--Kirillov dimension, which measures the growth of spaces
of invariants under small open subgroups, plays a key role in the study
of these representations.
 
In the unramified case, Breuil--Herzig--Hu--Morra--Schraen
\cite{BHH23,BHH25} introduced a graded ideal $J$ whose annihilation
condition yields the upper bound $n$ for Gelfand--Kirillov dimension.
Combined with a lower bound from patching, this gives dimension $n$
for representations arising from the mod-$p$ cohomology of
Shimura curves. For arbitrary $K$, Sorgdrager \cite{Sor26b} introduced
Casimir ideals and proved the same upper bound for
admissible Banach representations of $\GL_2(K)$ and $D^\times$
whose locally $\Q_p$-analytic vectors have an infinitesimal character.

On the lower-bound side, Hu--Morra--Schraen \cite{HMS14} proved that
the Iwasawa algebra of the upper unipotent subgroup of
$\GL_2(\cO_K)$ acts faithfully on every infinite-dimensional
irreducible smooth representation of $\GL_2(K)$. This suggests a
lower bound of $n$ for infinite-dimensional admissible smooth
$\GL_2(K)$-representations with a central character. Timmins
\cite{Tim26} proved that, when $K\ne\Q_p$, canonical dimension $1$
cannot occur for modules over the Iwasawa algebras of compact open
subgroups of $\SL_2(K)$. Thus the suggested lower bound holds
when $n=2$.

In this article, we prove that every infinite-dimensional finitely generated
module satisfying the Casimir annihilation condition has
Gelfand--Kirillov dimension $n$. 

Write $G=\GL_2(K)$ or $D^\times$, and fix a uniformizer $\varpi_K$
of $K$. We use the small uniform subgroup $H$ of the determinant-one
or reduced-norm-one group defined in Section~\ref{subsec:groups}.
Put $\Lambda=\Lambda_H:=\F[[H]]$. The weighted filtration on $\Lambda$
has associated graded algebra
\[
 B=B_H:=\gr_\omega\Lambda
   =\F[e_{k,j},f_{k,j},h_{k,j}:0\leq k\leq f-1,\ 0\leq j\leq e-1].
\]
The Casimir ideal $J\subset B$ is generated, for all $k$, by
the coefficients of
\[
 H_k(X)^2+4X E_k(X)F_k(X),
\]
where these polynomials are specified in
Section~\ref{sec:background}.
For a finitely generated $\Lambda$-module $M$, the condition
we consider is
\begin{equation}\label{eq:intro-casimir}
 J^N\gr_\omega M=0
 \qquad\text{for some }N\geq1.
\end{equation}
Here $\gr_\omega M$ is taken with respect to the natural good filtration
$F_\omega^nM=(F_\omega^n\Lambda)M$; see
Section~\ref{subsec:filtrations-dimension}. By
Lemma~\ref{lem:good-filtration-annihilation}, 
this is equivalent to annihilation by some power of $J$ for any
complete separated good module filtration compatible with
$F_\omega^\bullet\Lambda$. The exponent may depend on the filtration.

In practice, suppose that $K/\Q_p$ is unramified and $p\geq5$,
and take $C=1$ in the definition of $\omega$. For $M=\pi^\vee$,
where $\pi$ is an admissible smooth $G$-representation with a
central character, \eqref{eq:intro-casimir} is equivalent to
$\pi$ belonging to the category $\mathcal C$ defined by the graded
annihilation condition in \cite[\S~1]{Sor26a}; see
Section~\ref{subsec:casimir-ideal}. In the split case, this is a full
subcategory of the category denoted $\mathcal C$ in \cite{BHH25};
equality of the two categories is not known \cite[Warning~1]{Sor26a}.
For arbitrary $K$, the condition is also satisfied by
$M=(\mathcal V^\circ/\varpi_E\mathcal V^\circ)^\vee$, where
$\mathcal V$ is an infinite-dimensional admissible Banach
representation of $G$ over a finite extension $E/\Q_p$ whose
locally $\Q_p$-analytic vectors have an infinitesimal character,
$\varpi_E$ is a uniformizer of $E$, and $\mathcal V^\circ$ is any
$\widetilde H$-stable open bounded lattice; see
Corollary~\ref{cor:banach-gk-dimension}.

Now we state our main theorem. Choose a uniformizer $\varpi_D$ of $D$ with
$\varpi_D^2=\varpi_K$, and put
\[
 \Pi=\begin{pmatrix}0&1\\\varpi_K&0\end{pmatrix},\qquad
 \Theta=\begin{cases}
   \Pi,&G=\GL_2(K),\\
   \varpi_D,&G=D^\times,
 \end{cases}.
\]
Since $\Theta^2=\varpi_K$, the field $K'=K(\Theta)$ is a ramified
quadratic extension of $K$. Let $T=\exp(p^2\cO_K\Theta)$. This is a
torsion-free open subgroup of the norm-one principal units of
$K'$, contained in $H$, and $T\cong\Z_p^n$.
The associated norm-one algebraic $K$-torus has minimal splitting
field $K'$, so it is a ramified torus. By a slight abuse of terminology,
we also call $T$ the
ramified torus. 
Its Iwasawa algebra $\Lambda_T:=\F[[T]]$ is a power series ring
in $n$ variables.
Write $\delta_H(M)$ for the canonical dimension of $M$,
defined in Section~\ref{subsec:filtrations-dimension}.
Our main result uses only the $\Lambda$-module structure.

\begin{theorem}[Theorem~\ref{thm:main-lambda}]
\label{thm:intro-main}
Let $M$ be an infinite-dimensional finitely generated
$\Lambda$-module satisfying \eqref{eq:intro-casimir}. Then $M$ is
a finitely generated faithful $\Lambda_T$-module, and
\[
 \delta_H(M)=n.
\]
\end{theorem}

For an admissible smooth $G$-representation $\pi$ with a central
character, the pseudocompact dual $\pi^\vee$ is a finitely generated
$\Lambda$-module. Applying the theorem to $M=\pi^\vee$
gives the following corollary.

\begin{corollary}[Corollary~\ref{cor:main-gk-dimension}]
\label{cor:intro-gk-dimension}
Let $\pi$ be an infinite-dimensional admissible smooth
representation of $G$ over $\F$ with a central character. Suppose
that $\pi^\vee$ satisfies \eqref{eq:intro-casimir}.
Then $\pi^\vee$ is a finitely generated faithful $\Lambda_T$-module,
$\pi|_T$ is admissible, and
\[
 \GKdim(\pi)=n.
\]
\end{corollary}

The proof has two steps. First, we prove finite generation by the
same argument for both groups. The Casimir quotient $R:=B/J$ is
finite free of rank $4^n$ over
\[
 \gr_\omega \Lambda_T=B_T:=\F[e_{k,j}+f_{k,j}:k,j].
\]
The Casimir annihilation condition \eqref{eq:intro-casimir} then implies that
$\gr_\omega M$ is finite over $B_T$. Lifting generators gives
finite generation over $\Lambda_T$, hence the upper bound
$\delta_H(M)\leq n$, and reduces equality to faithfulness over $\Lambda_T$.

Second, we prove faithfulness using Frobenius pairs and the control
of ideals developed in \cite{AWZ08,AW09}. Conjugation by the subgroup
$\Gamma$ defined in Section~\ref{subsec:casimir-tower} gives sources of
derivations along a tower of open uniform subgroups indexed by powers
of $\varpi_K$. The resulting
control of localized annihilator ideals excludes nonzero
$\Lambda_T$-annihilators after inverting the polynomials $\Delta$ defined
in Section~\ref{subsec:casimir-local-control}. In the quaternionic
case, the open sets $\Delta\ne0$ cover $\operatorname{V}(J)$ away from
the origin. In the split case, the remaining subspaces are
$f_{k,j}=h_{k,j}=0$ for all $k,j$ and
$e_{k,j}=h_{k,j}=0$ for all $k,j$. We treat these subspaces using
the fact that the localized Iwasawa algebras of the upper and lower
unipotent subgroups have no proper nonzero $\Gamma$-stable ideals.

Since $M$ is finite and faithful over $\Lambda_T$, its support is all of
$\Spec \Lambda_T$. The dimension comparison therefore gives
$\delta_H(M)=\Krdim_{\Lambda_T} M=\Krdim \Lambda_T=n$.

Section~\ref{sec:background} recalls the groups, their filtered
Iwasawa algebras, the Casimir ideal, and the control theorem.
Section~\ref{sec:ramified-torus} constructs the ramified torus and
proves finite generation and the dimension comparison.
Section~\ref{sec:faithfulness} proves faithfulness.
Section~\ref{sec:conclusion} gives the main theorem and its
applications to smooth and Banach representations. The relation
with the unramified graded conditions of \cite{BHH23,HW24,Sor26a}
is explained in Section~\ref{subsec:casimir-ideal}.
In Appendix~\ref{app:canonical-diagrams}, Yongquan Hu proves that,
when $K\ne\Q_p$ and $\F$ is algebraically closed, the canonical
diagram of an admissible irreducible supersingular representation
satisfying \eqref{eq:intro-casimir} is $(\pi,\pi,\mathrm{id})$.

\paragraph*{Acknowledgments.} 
I am grateful to my advisor, Yongquan Hu, for suggesting the topic, for many useful discussions, and for writing the appendix to this paper.

\paragraph*{AI Disclosure.}
The methods used in this paper were discovered by GPT-6 Astra,  
which also helped draft the manuscript. 
The author has verified the results and proofs and revised the manuscript
to improve its clarity and readability. The author takes full responsibility
for the mathematical content and correctness of this paper.

\section{Background}\label{sec:background}

We fix notation and recall the filtered-algebra results used below.
Let $p\geq3$, and let $K/\Q_p$ be a finite extension with ramification
index $e$, residue degree $f$, and $n=[K:\Q_p]=ef$. Fix a uniformizer
$\varpi_K$ of $K$, and put $k_K=\cO_K/\varpi_K\cO_K$ and $q=p^f$.
When $K/\Q_p$ is unramified, we may take $\varpi_K=p$.
Let $K_0$ be the maximal unramified subfield of $K$, and let $L/K$
be the unramified quadratic extension, with residue field $k_L$.
Write $\tau$ for the nontrivial element of $\operatorname{Gal}(L/K)$.
Let $\F$ be any field of characteristic $p$ admitting an embedding
$\sigma_0:k_L\hookrightarrow\F$. Fix such an embedding
and put $\sigma_k=\sigma_0\circ\operatorname{Frob}_p^k$ for
$k\in\Z/2f\Z$. We use the same notation for their restrictions to $k_K$.
Throughout, embedding indices are read modulo $f$ on $k_K$ and
modulo $2f$ on $k_L$; the same convention applies to the corresponding
coordinate indices.

\subsection{The groups and their filtered Iwasawa algebras}
\label{subsec:groups}\label{subsec:graded-algebra}

Let $D$ be the non-split quaternion algebra over $K$. Choose a
Teichm\"uller unit $\eta\in\cO_L^{\times}$ such that
$\tau(\eta)=-\eta$ and $L=K(\eta)$. We fix a uniformizer
$\varpi_D$ of $D$ and the presentation
\begin{equation}\label{eq:quaternion-presentation}
 D=L\oplus\varpi_D L,\qquad \varpi_D^2=\varpi_K,\qquad
 \varpi_D a=\tau(a)\varpi_D\quad(a\in L).
\end{equation}
Its maximal order is $\cO_D=\cO_L\oplus\varpi_D\cO_L$.
We write $\Nrd:D\to K$ and $\Trd:D\to K$ for the reduced norm and trace.
As in \cite[\S~6.2]{Sor26b}, we use common notation for the two cases:
\[
 G=\begin{cases}
     \GL_2(K),&\text{in the split case},\\
     D^\times,&\text{in the quaternionic case}.
   \end{cases}
\]
Let $Z$ be the center of $G$, identified with $K^\times$, and let
$\nu_G:G\to K^\times$ be the determinant or reduced norm, respectively.
In the split case put
\[
 E_+=\begin{pmatrix}0&1\\0&0\end{pmatrix},\qquad
 E_-=\begin{pmatrix}0&0\\1&0\end{pmatrix},\qquad
 H_{\mathrm{diag}}=\begin{pmatrix}1&0\\0&-1\end{pmatrix}.
\]
Thus $\Pi=E_++\varpi_K E_-$.
Following \cite[\S~6.2]{Sor26b}, first define
\[
 \widetilde H^{1/p}=
 \begin{cases}
  \begin{pmatrix}
   1+p\varpi_K\cO_K&p\cO_K\\
   p\varpi_K\cO_K&1+p\varpi_K\cO_K
  \end{pmatrix},&\text{split},\\[6pt]
  1+p\varpi_D\cO_D,&\text{quaternionic},
 \end{cases}
\]
and put
\[
 H^{1/p}=\widetilde H^{1/p}\cap\ker\nu_G.
\]
Thus $\widetilde H^{1/p}$ is the full compact open subgroup denoted
by $H^{1/p}$ in that reference; our $H^{1/p}$ is its determinant-one
or reduced-norm-one subgroup.

Normalize the valuation on $K$ and $L$ by $v_p(p)=1$.
In the quaternionic case its extension to $D$ is
\[
 v_p(x)=\frac12v_p(\Nrd(x))\quad(x\in D^\times),
 \qquad v_p(\varpi_D)=\frac1{2e}.
\]
Choose a rational number $C$ such that
\begin{equation}\label{eq:valuation-constant-range}
 1-\frac1{p-1}<C<1-\frac1{p-1}+\frac1{2e}.
\end{equation}
In the split case write $g=\bigl(\begin{smallmatrix}a&b\\c&d\end{smallmatrix}\bigr)$.
The explicit $p$-valuation of \cite[(13)--(14)]{Sor26b}, with
the subtraction constant $C$ fixed above, is
\begin{equation}\label{eq:group-valuation}
 \omega_{1/p}(g)=
 \begin{cases}
  \min\left\{v_p(a-1),\ v_p(b)+\frac1{2e},
             v_p(c)-\frac1{2e},\ v_p(d-1)\right\}-C,
       &\text{split},\\[4pt]
  v_p(g-1)-C,&\text{quaternionic}.
 \end{cases}
\end{equation}
In either case, $\omega_{1/p}$ is a strictly saturated
$p$-valuation on $\widetilde H^{1/p}$ and restricts to
$H^{1/p}$. Now define
\begin{equation}\label{eq:groups-and-valuation}
 \widetilde H=(\widetilde H^{1/p})^p,\qquad
 H=(H^{1/p})^p,\qquad
 \omega=\omega_{1/p}|_H-1.
\end{equation}
Here the superscript $p$ denotes the subgroup of $p$th powers.
The shift by $1$ gives $\omega(g^p)=\omega_{1/p}(g)$ for
$g\in H^{1/p}$. We use this shifted valuation $\omega$ for the
weighted filtration below. Passing to $p$th powers also makes the associated
graded group algebras commutative for all the radius filtrations
used in the Banach application; see \cite[Lemma~4.3 and \S~6.2]{Sor26b}.

To describe $\log H$, put
\[
 \mathfrak j=
 \begin{cases}
  \begin{pmatrix}
   \varpi_K\cO_K&\cO_K\\
   \varpi_K\cO_K&\varpi_K\cO_K
  \end{pmatrix},&\text{split},\\[6pt]
  \varpi_D\cO_D
    =\varpi_K\cO_L\oplus\varpi_D\cO_L,&\text{quaternionic}.
 \end{cases}
\]
Thus $\widetilde H^{1/p}=1+p\mathfrak j$.
For $m\geq1$, the exponential and logarithm give inverse bijections
between $p^m\mathfrak j$ and $1+p^m\mathfrak j$; see
\cite[\S\S~3.1 and 6.2]{Sor26b}.
Write $\operatorname{tr}_G$ for the ordinary or reduced trace,
respectively, and put $\mathfrak j_0=\mathfrak j\cap\ker\operatorname{tr}_G$.
Using $\nu_G(\exp X)=\exp(\operatorname{tr}_G X)$ and
$(\exp X)^p=\exp(pX)$, we obtain
$H^{1/p}=\exp(p\mathfrak j_0)$ and $\log H=\mathfrak l:=p^2\mathfrak j_0$.
Explicitly,
\begin{equation}\label{eq:compact-subgroups}
 \mathfrak l=
 \begin{cases}
  p^2(\cO_KE_+\oplus\varpi_K\cO_KE_-
                   \oplus\varpi_K\cO_KH_{\mathrm{diag}}),&\text{split},\\
  p^2(\cO_K\varpi_D\oplus\cO_K\eta\varpi_D
                   \oplus\varpi_K\cO_K\eta),&\text{quaternionic},
 \end{cases}
 \qquad H=\exp(\mathfrak l).
\end{equation}
In both cases $\mathfrak l$ is a finite free $\Z_p$-Lie lattice
of rank $3n$ with $[\mathfrak l,\mathfrak l]\subseteq p\mathfrak l$.
Hence $H$ is uniform by \cite[Thm.~9.10]{DDMS99}.
Put
\[
 Z_H=\exp(p^2\varpi_K\cO_K\cdot1).
\]
Then $\widetilde H$ is a compact open subgroup of $G$ of dimension $4n$,
and
\begin{equation}\label{eq:norm-one-realization}
 \widetilde H=H\times Z_H,\qquad H=\widetilde H\cap\ker\nu_G.
\end{equation}
Indeed, the logarithm of $\widetilde H$ is
$\mathfrak l\oplus p^2\varpi_K\cO_K\cdot1$.
The map $\nu_G$ is trivial on $H$ and is squaring on
$Z_H$; squaring is an automorphism of $Z_H$. These are the small
groups and their central quotients used in \cite[\S\S~6.2 and 6.5]{Sor26b}.

We give $\F$ the discrete topology. The Iwasawa algebra of a
profinite group $P_0$ over $\F$ is
\[
 \F[[P_0]]=\varprojlim_U\F[P_0/U],
\]
where $U$ runs over its open normal subgroups. We put
\[
 \Lambda=\Lambda_H:=\F[[H]],\qquad \fm=\ker(\Lambda\longrightarrow\F).
\]
The ring $\Lambda$ is noetherian and local, with maximal ideal $\fm$.
Its inverse-limit topology is the $\fm$-adic topology. It is a
pseudocompact $\F$-algebra: it is an inverse limit of
finite-dimensional discrete $\F$-algebras. These assertions hold
without a finiteness assumption on $\F$; see
\cite[Prop.~1.9 and Thm.~3.1]{Koh17}.

We equip every finitely generated $\Lambda$-module $M$ with its
$\fm$-adic topology. With this topology, $M$ is a complete separated
pseudocompact $\Lambda$-module; see \cite[\S~1]{Koh17}.
For a smooth representation $\pi$ over $\F$, write
$\pi^\vee=\Hom_{\F}(\pi,\F)$, with the topology of pointwise
convergence and the contragredient action. Its continuous
$\F$-linear dual is $\pi$. This duality identifies smooth
$H$-representations with pseudocompact $\Lambda$-modules, and $\pi$ is
admissible if and only if $\pi^\vee$ is finitely generated over
$\Lambda$; see \cite[Thm.~1.5, Cor.~1.8 and Prop.~1.9]{Koh17}.
If $\pi$ is a smooth $G$-representation with a central character, then
$Z_H$ acts trivially: the image of a smooth character of $Z_H$
is a finite $p$-group, whereas $\F^\times$ has no nontrivial
$p$-power torsion. If $\pi$ is also admissible, then $\pi|_H$ is admissible:
for every open subgroup $U\subseteq H$, the subgroup $UZ_H$ is open
in $G$ and $\pi^U=\pi^{UZ_H}$. Thus $\pi^\vee$ is finite over
$\Lambda$.

We now describe the weighted filtration on $\Lambda$ induced by the
$p$-valuation $\omega$ in \eqref{eq:groups-and-valuation} and its
associated graded algebra.
Choose $\alpha\in\cO_{K_0}^{\times}$ such that
$1,\alpha,\ldots,\alpha^{f-1}$ is a $\Z_p$-basis of $\cO_{K_0}$.
For $0\leq i\leq f-1$ and $0\leq j\leq e-1$, define the split group elements
\begin{align*}
 u^+_{i,j}&=\exp(p^2\alpha^i\varpi_K^jE_+)
    =\begin{pmatrix}1&p^2\alpha^i\varpi_K^j\\0&1\end{pmatrix},\\[4pt]
 u^-_{i,j}&=\exp(p^2\alpha^i\varpi_K^{j+1}E_-)
    =\begin{pmatrix}1&0\\p^2\alpha^i\varpi_K^{j+1}&1\end{pmatrix},\\[4pt]
 u^0_{i,j}&=\exp(p^2\alpha^i\varpi_K^{j+1}H_{\mathrm{diag}})\\
    &=\begin{pmatrix}
       \exp(p^2\alpha^i\varpi_K^{j+1})&0\\
       0&\exp(-p^2\alpha^i\varpi_K^{j+1})
      \end{pmatrix}.
\end{align*}
The elements $u^+_{i,j},u^-_{i,j},u^0_{i,j}$ give an ordered
basis of $H$. In the quaternionic case, use
the ordered basis
\begin{align*}
 u^+_{i,j}&=\exp(p^2\alpha^i\varpi_K^j\varpi_D),\\
 u^-_{i,j}&=\exp(p^2\alpha^i\varpi_K^j\eta\varpi_D),\\
 u^0_{i,j}&=\exp(p^2\alpha^i\varpi_K^{j+1}\eta).
\end{align*}
In both cases the weights of these group elements are
\[
 \begin{aligned}
 \omega(u^+_{i,j})=\omega(u^-_{i,j})&=c+j/e,\\
 \omega(u^0_{i,j})&=c+(j+1/2)/e,
 \end{aligned}
\]
where
\[
 c=1+\frac1{2e}-C.
\]
By \eqref{eq:valuation-constant-range}, all these weights lie strictly
between $1/(p-1)$ and $p/(p-1)$;
the largest is $c+(e-1/2)/e$. Order the elements
$u^+_{i,j},u^-_{i,j},u^0_{i,j}$ as
$g_1,\ldots,g_{3n}$. Every element of $H$ has a unique expression
$g_1^{z_1}\cdots g_{3n}^{z_{3n}}$ with $z_i\in\Z_p$.
Put $b_i=g_i-1\in\Lambda$. Every element of $\Lambda$ has a
unique ordered power series expansion
\[
 x=\sum_{\mathbf r\in\N^{3n}}c_{\mathbf r}b^{\mathbf r},
 \qquad c_{\mathbf r}\in\F,\qquad
 b^{\mathbf r}=b_1^{r_1}\cdots b_{3n}^{r_{3n}}.
\]
Fix a positive integer $\ell_\omega$ clearing the denominators of all
basis weights.
Set
\[
 w_i=\ell_\omega\,\omega(g_i),\qquad
 a_j=\ell_\omega\left(c+\frac je\right),\qquad
 \nu=\frac{\ell_\omega}{2e}.
\]
The numbers $w_i,a_j,\nu$ are positive integers. The integer
weights are $a_j$ on $u^\pm_{i,j}-1$ and $a_j+\nu$ on
$u^0_{i,j}-1$ in both cases.
For $\mathbf r=(r_1,\ldots,r_{3n})$, define its monomial weight by
\[
 \tau_{\mathbf r}:=\sum_{i=1}^{3n}w_i r_i
     =\ell_\omega\sum_{i=1}^{3n}r_i\omega(g_i).
\]
The weighted filtration on $\Lambda$ is, for $m\in\Z$,
\begin{equation}\label{eq:weighted-ring-filtration}
 F_\omega^m\Lambda=
 \left\{\sum_{\mathbf r\in\N^{3n}}c_{\mathbf r}b^{\mathbf r}:
 c_{\mathbf r}=0\ \text{if }\tau_{\mathbf r}<m\right\}.
\end{equation}
Thus $F_\omega^m\Lambda$ consists of the series all of whose nonzero
monomials have weight at least $m$; in particular,
$F_\omega^0\Lambda=\Lambda$ and $F_\omega^1\Lambda=\fm$.
This is a decreasing multiplicative filtration. Since there are
finitely many positive generator weights, it defines the same
topology as the $\fm$-adic filtration and is complete and separated.
We use the integer rescaling fixed above throughout, and write
\[
 B=B_H:=\gr_\omega\Lambda
   =\bigoplus_{m\geq0}F_\omega^m\Lambda/F_\omega^{m+1}\Lambda,
\]
with weighted degree $\deg_\omega=m$ on
$F_\omega^m\Lambda/F_\omega^{m+1}\Lambda$. Thus $\omega$ denotes
the group valuation, while $F_\omega^\bullet\Lambda$ is its
integer-rescaled ring filtration. We omit the subscript $H$ in
$\Lambda_H$ and $B_H$; for the other subgroup algebras below,
matching subscripts on $\Lambda$ and $B$ indicate the filtered
algebra and its associated graded algebra, respectively.

We now define the graded coordinates explicitly.
Let $\beta_0,\ldots,\beta_{f-1}\in\cO_{K_0}$ be the basis
trace-dual to $1,\alpha,\ldots,\alpha^{f-1}$, so that
\[
 \operatorname{Tr}_{K_0/\Q_p}(\beta_i\alpha^{i'})=\delta_{i,i'}
 \qquad(0\leq i,i'\leq f-1).
\]
This integral dual basis exists because $K_0/\Q_p$ is unramified.
Bars denote reduction to the residue field. The matrix
$(\sigma_k(\overline\beta_i))_{k,i}$ satisfies
\[
 \sum_{i=0}^{f-1}\sigma_k(\overline\beta_i)
                     \sigma_{k'}(\overline\alpha)^i
     =\delta_{k,k'}\qquad(0\leq k,k'\leq f-1),
\]
so it is the inverse of $(\sigma_{k'}(\overline\alpha)^i)_{i,k'}$.
We use this matrix to change the basis of principal symbols below.
For $b\in F_\omega^n\Lambda\setminus F_\omega^{n+1}\Lambda$, write
$\gr_\omega(b)=b+F_\omega^{n+1}\Lambda$ for its principal symbol
in the degree-$n$ component of $B$. Thus $\gr_\omega(b)$ denotes
the principal symbol of an element, whereas $\gr_\omega\Lambda$
and $\gr_\omega M$ denote associated graded objects.

In the split case, for $0\leq k\leq f-1$ and $0\leq j\leq e-1$, put
\begin{equation}\label{eq:split-graded-coordinates}
 \begin{aligned}
 e_{k,j}&=\sum_{i=0}^{f-1}\sigma_k(\overline\beta_i)
                            \gr_\omega(u^+_{i,j}-1),\\
 f_{k,j}&=\sum_{i=0}^{f-1}\sigma_k(\overline\beta_i)
                            \gr_\omega(u^-_{i,j}-1),\\
 h_{k,j}&=\sum_{i=0}^{f-1}\sigma_k(\overline\beta_i)
                            \gr_\omega(u^0_{i,j}-1).
 \end{aligned}
\end{equation}

In the quaternionic case, for $0\leq k\leq 2f-1$ and $0\leq j\leq e-1$, put
\[
 x_{k,j}=\frac12\sum_{i=0}^{f-1}\left(
    \sigma_k(\overline\beta_i)\gr_\omega(u^+_{i,j}-1)
   +\sigma_k(\overline\beta_i/\overline\eta)
                           \gr_\omega(u^-_{i,j}-1)\right).
\]
For $0\leq k\leq f-1$ define
\begin{equation}\label{eq:quaternionic-graded-coordinates}
 \begin{aligned}
 e_{k,j}&=x_{k,j},\qquad f_{k,j}=x_{k+f,j},\\
 h_{k,j}&=\sum_{i=0}^{f-1}
             \sigma_k(\overline\beta_i/\overline\eta)
                           \gr_\omega(u^0_{i,j}-1).
 \end{aligned}
\end{equation}
In fact, $e_{k,j},f_{k,j},h_{k,j}$ are obtained from the formulas in
\cite[Lemma~6.3]{Sor26b} by setting $N=N'=0$ and reducing modulo
the principal symbol of the coefficient-field uniformizer.

\begin{proposition}\label{prop:graded-algebra}
In both cases there is an isomorphism of graded $\F$-algebras
\begin{equation}\label{eq:graded-presentation}
 B=B_H:=\gr_\omega\Lambda
   =\F[e_{k,j},f_{k,j},h_{k,j}:0\leq k\leq f-1,\ 0\leq j\leq e-1].
\end{equation}
With $a_0,\ldots,a_{e-1},\nu$ as defined above,
\[
 \deg_\omega e_{k,j}=\deg_\omega f_{k,j}=a_j,\qquad
 \deg_\omega h_{k,j}=a_j+\nu.
\]
\end{proposition}

\begin{proof}
The polynomial presentation follows from \cite[Lemma~6.3]{Sor26b}
with $N=N'=0$ and the variables $z_{k,j}$ for the central subgroup
$Z_H$ omitted.
Since the coordinate formulas are defined over $\sigma_0(k_L)$,
the presentation descends to this field and then extends to $\F$.
The degrees follow from the weights of the ordered basis elements.
\end{proof}

The central subgroup $Z_H$ has ordered basis
$g^Z_{i,j}:=\exp(p^2\alpha^i\varpi_K^{j+1})$
for $0\leq i\leq f-1$ and $0\leq j\leq e-1$.
Put $\widetilde\Lambda=\F[[\widetilde H]]$. Extend the weighted
filtration from $\Lambda$ to $\widetilde\Lambda$ by assigning
weight $a_j+\nu$ to
$g^Z_{i,j}-1$. Writing $\gr_\omega$ also for principal symbols
in this extended filtration, define, for $0\leq k\leq f-1$ and $0\leq j\leq e-1$,
\[
 z_{k,j}=\sum_{i=0}^{f-1}\sigma_k(\overline\beta_i)
    \gr_\omega(g^Z_{i,j}-1).
\]
This extended filtration on $\widetilde\Lambda$ gives
\[
 \widetilde B:=\gr_\omega\widetilde\Lambda
    =B[z_{k,j}:0\leq k\leq f-1,\ 0\leq j\leq e-1],\qquad
 \deg_\omega z_{k,j}=a_j+\nu.
\]

\subsection{The Casimir ideal}\label{subsec:casimir-ideal}

Let $g(X)\in\cO_{K_0}[X]$ be the Eisenstein minimal polynomial of
$\varpi_K$. Define $\gamma_j\in K$ by
\[
 \frac{g(X)}{(X-\varpi_K)g'(\varpi_K)}
       =\sum_{j=0}^{e-1}\gamma_jX^j.
\]
Put $R_K=e\,v_p(g'(\varpi_K))+1-e$. By
\cite[Lemma~6.2]{Sor26b}, $R_K\in\Z$ and
$v_p(\gamma_j)=-(R_K+j)/e$.
Put
\[
 \mu_{k,j}=\sigma_k\bigl(\overline{\varpi_K^{R_K+j}\gamma_j}\bigr)
          \in\F^\times.
\]
Since $\varpi_K^{R_K+j}\gamma_j\in\cO_K^\times$, we have
$\mu_{k,j}\in\sigma_0(k_K)^\times$. In particular,
$\mu_{k,j}^{p^f}=\mu_{k,j}$. Define
\begin{align*}
 E_k(X)&=\sum_{j=0}^{e-1}\mu_{k,j}e_{k,j}X^{e-1-j},\\
 F_k(X)&=\sum_{j=0}^{e-1}\mu_{k,j}f_{k,j}X^{e-1-j},\\
 H_k(X)&=\sum_{j=0}^{e-1}\mu_{k,j}h_{k,j}X^{e-1-j}.
\end{align*}
The Casimir ideal is
\begin{equation}\label{eq:general-casimir-ideal}
 J=
 \bigl([X^m](H_k(X)^2+4X E_k(X)F_k(X)):
           0\leq k\leq f-1,\ 0\leq m\leq 2e-1\bigr)\subset B.
\end{equation}
Here $[X^m]$ means the coefficient of $X^m$; all $2e$
coefficients are included for each $k$. The ideal $J$ is homogeneous
for $\deg_\omega$: the generator
$[X^m](H_k(X)^2+4X E_k(X)F_k(X))$ has weighted degree
$2a_0+2\nu(2e-1-m)$.
For the auxiliary grading defined by
$\deg_{\mathrm{tot}}e_{k,j}=\deg_{\mathrm{tot}}f_{k,j}=\deg_{\mathrm{tot}}h_{k,j}=1$,
every generator of $J$ is homogeneous of degree $2$.
We write $R=B/J$ for the Casimir quotient.

Our definition of $J$ agrees with Sorgdrager's Casimir ideal in
\cite[proof of Thm.~6.2, (22)--(24)]{Sor26b} after quotienting by
the central variables $z_{k,j}$. In the coordinates fixed above,
reducing the coefficients of the normalized Casimir symbols in
\cite[(22)]{Sor26b} at $N=0$ gives
\[
 [X^m]\left(\frac12 H_k(X)^2+2X E_k(X)F_k(X)\right),
 \qquad 0\leq k\leq f-1,\quad 0\leq m\leq 2e-1.
\]
These are one half of the generators in
\eqref{eq:general-casimir-ideal}, so they generate $J$ since
$p\ne2$. Sorgdrager's ideal for the full group $\widetilde H$
is therefore
\[
 J\,\widetilde B+(z_{k,j}:0\leq k\leq f-1,\ 0\leq j\leq e-1).
\]

We now compare the Casimir ideal with the unramified graded ideal
of \cite{BHH23} and its quaternionic analogue in \cite{HW24}.
Suppose that $p\geq5$, that $K/\Q_p$ is unramified, and that
$\varpi_K=p$. For this comparison we take $C=1$, which lies in
\eqref{eq:valuation-constant-range} since $e=1$ and $p\geq5$.
Then $c=1/2$, so the unrescaled weights of
$e_{k,0},f_{k,0},h_{k,0}$ are $1/2,1/2,1$, respectively.
Also $n=f$ and
$g(X)=X-p$, so $\gamma_0=1$, $R_K=0$, and $\mu_{k,0}=1$.
Consequently \eqref{eq:general-casimir-ideal} becomes
\[
 J=(h_{k,0}^2,e_{k,0}f_{k,0}:0\leq k\leq f-1),\qquad
 \sqrt J=(h_{k,0},e_{k,0}f_{k,0}:0\leq k\leq f-1).
\]

Let $I_1$ be the upper pro-$p$ Iwahori subgroup of $\GL_2(K)$,
put $U_D^1=1+\varpi_D\cO_D$ and $Z_1=1+p\cO_K$, and write
\[
 G_u=\begin{cases}I_1/Z_1,&\text{split},\\
                   U_D^1/Z_1,&\text{quaternionic},\end{cases}
 \qquad \Lambda_u=\F[[G_u]],\qquad
 \mathfrak m_u=\ker(\Lambda_u\to\F).
\]
There are isomorphisms
\[
 G_u\cong
 \begin{cases}
  I_1\cap\mathrm{SL}_2(K),&\text{split},\\
  \{g\in U_D^1:\Nrd(g)=1\},&\text{quaternionic}.
 \end{cases}
\]
sending $gZ_1$ to $\nu_G(g)^{-1/2}g$. 
Under this identification, $H=G_u^{p^2}$; see
\cite[Prop.~2]{Sor26a}.
Write $\gr_{\mathfrak m_u}$ for the principal-symbol map for the
$\mathfrak m_u$-adic filtration, and write $[\lambda]$ for the
Teichm\"uller lift of a residue-field element $\lambda$.
In the split case, use the normalization in
\cite[(3.1), (3.17)]{BHH25}: for $0\leq k\leq f-1$, put
\[
 \begin{aligned}
 Y_k&=\sum_{\lambda\in k_K^\times}\sigma_k(\lambda)^{-1}
             \begin{pmatrix}1&[\lambda]\\0&1\end{pmatrix},\\
 Z_k&=\sum_{\lambda\in k_K^\times}\sigma_k(\lambda)^{-1}
             \begin{pmatrix}1&0\\p[\lambda]&1\end{pmatrix},\\
 y_k&=\gr_{\mathfrak m_u}(Y_k),\qquad z_k=\gr_{\mathfrak m_u}(Z_k),\qquad
 h_k=[y_k,z_k].
 \end{aligned}
\]
In the quaternionic case, follow \cite[\S~2.3]{HW24}: for
$0\leq k\leq 2f-1$, put
\[
 Y_k=\sum_{\lambda\in k_L^\times}\sigma_k(\lambda)^{-1}
                    (1+\varpi_D[\lambda]),
 \qquad y_k=\gr_{\mathfrak m_u}(Y_k).
\]
For $0\leq k\leq f-1$, set
\[
 z_k=y_{k+f},\qquad h_k=[y_k,z_k].
\]
Our $h_k$ is the generator denoted $h'_k$ in
\cite[Prop.~2.9]{HW24}.
The algebra $\gr_{\mathfrak m_u}\Lambda_u$ is generated by
$y_k,z_k$ of degree $1$ and central $h_k$ of degree $2$;
all brackets between degree-one generators are central.
The two-sided ideal used in \cite[\S~5.3]{BHH23} and
\cite[Cor.~2.10]{HW24} is
\[
 J_u=(h_k,y_kz_k:0\leq k\leq f-1)\subset\gr_{\mathfrak m_u}\Lambda_u.
\]

Taking $C=1$, the $\omega$-filtration on $\Lambda$ agrees,
up to a common rescaling, with the filtration induced from the
$\mathfrak m_u$-adic filtration on $\Lambda_u$;
see \cite[Prop.~5.3.3]{BHH23} and \cite[Lemma~2.4]{HW24}.
Rescale the degrees on $B$ so that
$e_{k,0},f_{k,0}$ have degree $p^2$ and $h_{k,0}$ has degree $2p^2$.
We identify $B$ with its image in
$\gr_{\mathfrak m_u}\Lambda_u$.
In the quaternionic case, extend $z_k,h_k$ by
$z_{k+f}=y_k$ and $h_{k+f}=-h_k$.
We now compute the $p^2$th powers of these generators in terms
of the coordinates defined in \eqref{eq:split-graded-coordinates}
and \eqref{eq:quaternionic-graded-coordinates}.
In the split case put $g_i^+=\exp(\alpha^iE_+)$.
The reduction of $[\lambda]$ in the basis
$1,\overline\alpha,\ldots,\overline\alpha^{f-1}$ gives
\[
 \gr_{\mathfrak m_u}\bigl(\exp([\lambda]E_+)-1\bigr)
 =\sum_{i=0}^{f-1}
   \operatorname{Tr}_{k_K/\mathbf F_p}(\overline\beta_i\lambda)
                         \gr_{\mathfrak m_u}(g_i^+-1).
\]
For $k'=k$, the inner sum below equals $q-1=-1$ in $\F$;
for $k'\ne k$, it is zero. Expanding the trace therefore gives
\[
 \begin{aligned}
 &\sum_{\lambda\in k_K^\times}\sigma_k(\lambda)^{-1}
              \operatorname{Tr}_{k_K/\mathbf F_p}(\overline\beta_i\lambda)\\
 &\qquad=\sum_{k'=0}^{f-1}\sigma_{k'}(\overline\beta_i)
             \sum_{\lambda\in k_K^\times}
                  \frac{\sigma_{k'}(\lambda)}{\sigma_k(\lambda)}
       =-\sigma_k(\overline\beta_i).
 \end{aligned}
\]
Consequently
$y_k=-\sum_i\sigma_k(\overline\beta_i)\gr_{\mathfrak m_u}(g_i^+-1)$.
The elements $g_i^+$ commute, $(g_i^+)^{p^2}=u^+_{i,0}$,
and $g^{p^2}-1=(g-1)^{p^2}$ in characteristic $p$.
Since $\sigma_{k-2}(b)^{p^2}=\sigma_k(b)$ for $b\in k_L$, it follows that
\[
 \begin{aligned}
 y_{k-2}^{p^2}
 &=-\sum_{i=0}^{f-1}\sigma_{k-2}(\overline\beta_i)^{p^2}
                                  \gr_{\mathfrak m_u}(g_i^+-1)^{p^2}\\
 &=-\sum_{i=0}^{f-1}\sigma_k(\overline\beta_i)
                                  \gr_{\mathfrak m_u}(u^+_{i,0}-1)
   =-e_{k,0}.
 \end{aligned}
\]
Replacing $E_+$ by $pE_-$ in the same calculation gives
$z_{k-2}^{p^2}=-f_{k,0}$.
The relation $[E_+,pE_-]=pH_{\mathrm{diag}}$ and the Lazard
isomorphism of \cite[(43)]{BHH23} give
\[
 h_k=[y_k,z_k]
     =\sum_{i=0}^{f-1}\sigma_k(\overline\beta_i)
              \gr_{\mathfrak m_u}\bigl(\exp(p\alpha^iH_{\mathrm{diag}})-1\bigr),
 \qquad h_{k-2}^{p^2}=h_{k,0}.
\]
In the quaternionic case put
\[
 g_i^1=\exp(\alpha^i\varpi_D),\qquad
 g_i^2=\exp(\alpha^i\eta\varpi_D),\qquad
 g_i^0=\exp(p\alpha^i\eta).
\]
The norm-one representative of $1+\varpi_D[\lambda]$
has logarithm congruent to
$\varpi_D[\lambda]=[\lambda^q]\varpi_D$ modulo
$\varpi_D^3\cO_D$. Since the trace-dual basis to
$\overline\alpha^i,\overline\alpha^i\overline\eta$ in $k_L$
is $\overline\beta_i/2,\overline\beta_i/(2\overline\eta)$,
the same finite-field sum calculation gives
\[
 \begin{aligned}
 y_k&=-\frac12\sum_{i=0}^{f-1}\Bigl(
       \sigma_k(\overline\beta_i)\gr_{\mathfrak m_u}(g_i^1-1)
       -\sigma_k(\overline\beta_i/\overline\eta)
                                  \gr_{\mathfrak m_u}(g_i^2-1)\Bigr),\\
 z_k&=-\frac12\sum_{i=0}^{f-1}\Bigl(
       \sigma_k(\overline\beta_i)\gr_{\mathfrak m_u}(g_i^1-1)
       +\sigma_k(\overline\beta_i/\overline\eta)
                                  \gr_{\mathfrak m_u}(g_i^2-1)\Bigr).
 \end{aligned}
\]
The relation $[\varpi_D,\eta\varpi_D]=-2p\eta$ then gives
\[
 h_k=[y_k,z_k]
     =-\sum_{i=0}^{f-1}\sigma_k(\overline\beta_i/\overline\eta)
                                      \gr_{\mathfrak m_u}(g_i^0-1).
\]
Their group powers satisfy
\[
 (g_i^1)^{p^2}=u^+_{i,0},\qquad
 (g_i^2)^{p^2}=u^-_{i,0},\qquad
 (g_i^0)^{p^2}=u^0_{i,0}.
\]
The $p$th-power map is
additive on the span of the degree-one generators since their
brackets are central. Applying this and
$\sigma_{k-2}(b)^{p^2}=\sigma_k(b)$ to the formulas above gives
\[
 \begin{aligned}
 y_{k-2}^{p^2}
 &=-\frac12\sum_{i=0}^{f-1}\Bigl(
       \sigma_k(\overline\beta_i)\gr_{\mathfrak m_u}(u^+_{i,0}-1)
       -\sigma_k(\overline\beta_i/\overline\eta)
                                  \gr_{\mathfrak m_u}(u^-_{i,0}-1)\Bigr)\\
 &=-x_{k+f,0}=-f_{k,0}.
 \end{aligned}
\]
The same calculation yields
\[
 z_{k-2}^{p^2}=-e_{k,0},\qquad h_{k-2}^{p^2}=-h_{k,0}.
\]
With the definitions in \cite[\S~2.3]{HW24}, the powers
$y_{k-2}^{p^2}$ and $z_{k-2}^{p^2}$ correspond to $-f_{k,0}$ and
$-e_{k,0}$. With the definitions in \cite[(3.1), (3.17)]{BHH25},
they correspond to $-e_{k,0}$ and $-f_{k,0}$.
These sign changes and the exchange of $e_{k,0}$ with $f_{k,0}$
leave the ideal $(h_{k,0},e_{k,0}f_{k,0}:k)$ unchanged.
In particular, we obtain
\[
 \sqrt J=(h_k^{p^2},y_k^{p^2}z_k^{p^2}:0\leq k\leq f-1)
                         \subset B.
\]
This is the ideal denoted $J_N$ in \cite[\S~5]{Sor26a} for $N=2$.
Since $J\subseteq\sqrt J$
and $(\sqrt J)^{f+1}\subseteq J$, annihilation by a power
of either ideal is equivalent.
The comparison in
\cite[\S~6.4]{Sor26b}, based on
\cite[\S~6, Lemmas~1--4]{Sor26a}, gives, for every finitely
generated $\Lambda_u$-module $M$,
\begin{equation}\label{eq:bhh-casimir-annihilation}
 \begin{aligned}
 &J_u^a\gr_{\mathfrak m_u}M=0\quad\text{for some }a\geq1\\
 &\qquad\Longleftrightarrow\quad
 J^b\gr_\omega(M|_H)=0
                       \quad\text{for some }b\geq1.
 \end{aligned}
\end{equation}

\subsection{Why we do not use \texorpdfstring{$\mathfrak m$}{m}-adic filtration}
\label{subsec:why-not-augmentation}

We explain in this subsection that for arbitrary extensions, it is more natural to use the uniform group $H$ and its weighted
filtration, differing from \cite{BHH23} in the unramified case. First, for ramified $K$, the upper pro-$p$ Iwahori subgroup
$I_1\subset\GL_2(K)$, modulo its pro-$p$ centre $Z_1=I_1\cap Z$,
need not be $p$-torsion-free. For example, for any $p\geq3$,
let $K=\Q_p(\zeta_p+\zeta_p^{-1})$ and $t=\zeta_p+\zeta_p^{-1}$,
where $\zeta_p$ is a primitive $p$th root of unity. Then
\[
 g=\begin{pmatrix}1&-1\\2-t&t-1\end{pmatrix}\in I_1\cap\SL_2(K)
\]
has order $p$. Thus it is natural for extensions with arbitrary
ramification to use a small uniform subgroup, such as
$H$ in \eqref{eq:compact-subgroups}.

Second, even for $H$, the $\fm$-adic filtration does not make the
required Casimir annihilation condition hold even for principal series.
To see this, take $K=\Q_p(\varpi_K)$ with $\varpi_K^2=p$.
Write $e_j,f_j,h_j$ ($j=0,1$) for the augmentation principal symbols
of $u^+_{0,j}-1,u^-_{0,j}-1,u^0_{0,j}-1$, respectively. Uniformity gives
\[
 \gr_{\fm}\Lambda=\F[e_0,e_1,f_0,f_1,h_0,h_1],
 \qquad \deg_{\fm}e_j=\deg_{\fm}f_j=\deg_{\fm}h_j=1.
\]
Thus $\gr_\omega\Lambda$ and $\gr_{\fm}\Lambda$ have the same
underlying polynomial algebra in the chosen coordinates, equipped
with the weighted grading and the standard grading, respectively.
The coefficient formula \eqref{eq:general-casimir-ideal} also defines
a Casimir ideal in $\gr_{\fm}\Lambda$, which we again denote by $J$.
The two Casimir ideals coincide as ideals of this polynomial algebra;
only their gradings differ.

Consider the smooth mod-$p$ principal series
\[
 \pi=\operatorname{Ind}_B^{\GL_2(K)}\mathbf1
     \cong C^\infty(\mathbf P^1(K),\F),
 \qquad M=\pi^\vee.
\]
Here $B$ is the upper triangular Borel subgroup and induction is
unnormalized. The representation $\pi$ is infinite-dimensional and
admissible, with trivial central character. 
It is directly checked that the orbit of $x=[1:1]$ is 
$\mathcal O=1+p^2\cO_K\subset \mathbf P^1(K)$. This orbit is
clopen, so, with $H_x=\operatorname{Stab}_H(x)$,
\[
 Q:=C^\infty(\mathcal O,\F)^\vee
   \cong\Lambda/\overline{\Lambda\fm_{H_x}}
\]
is a direct summand of $M$ as a $\Lambda$-module; here $\fm_{H_x}$
is the augmentation ideal of $\F[[H_x]]$.

The point $x=[1:1]$ corresponds to the line $K\cdot(1,1)\subset K^2$.
Thus $h\in H$ belongs to $H_x$ exactly when $h$ maps this line to
itself. The convergent
logarithm and exponential series show that $h$ preserves this line
if and only if $\log h$ does. Hence $\log H_x$ consists precisely of
the elements of $\mathfrak l$ preserving $K\cdot(1,1)$. Write such an element as
\[
 X=p^2\begin{pmatrix}\varpi_K a&b\\
                      \varpi_K c&-\varpi_K a\end{pmatrix},
 \qquad a,b,c\in\cO_K.
\]
The condition is that the two coordinates of $X(1,1)$ coincide:
\[
 X\binom11=p^2\binom{\varpi_K a+b}{\varpi_K c-\varpi_K a}
       \in K\cdot\binom11
 \quad\Longleftrightarrow\quad b=\varpi_K(c-2a).
\]
Substituting this expression for $b$ gives
\[
 \log H_x=p^2\varpi_K\bigl(
    \cO_K(H_{\mathrm{diag}}-2E_+)
       \oplus\cO_K(E_++E_-)\bigr).
\]
A $\Z_p$-basis of $\log H_x$ and their images in
$\fm/\fm^2$ is given by 
\[
 \begin{aligned}
 v_1&=p^2\varpi_K(H_{\mathrm{diag}}-2E_+)
       &&\longmapsto h_0-2e_1,\\
 v_2&=\varpi_Kv_1=p^3H_{\mathrm{diag}}-2p^3E_+
       &&\longmapsto h_1,\\
 v_3&=p^2\varpi_K(E_++E_-)
       &&\longmapsto f_0+e_1,\\
 v_4&=\varpi_Kv_3=p^3E_++p^3E_-
       &&\longmapsto f_1.
 \end{aligned}
\]
In the second and fourth rows, $p^3E_+=p(p^2E_+)\in p\mathfrak l$,
contributing zero to the degree-one symbol.
Therefore we have 
\[
 \gr_{\fm}Q\cong
 \frac{\F[e_0,e_1,f_0,f_1,h_0,h_1]}
      {(h_0-2e_1,h_1,f_0+e_1,f_1)}
 \cong\F[e_0,e_1]
\]
as $\gr_{\fm}\Lambda$-modules. 
We may compute in our case that $g(X)=X^2-p$ 
and $\mu_{0,0}=\mu_{0,1}=1/2$. Consequently,
$e_0f_0=[X^3](H_0(X)^2+4XE_0(X)F_0(X))$ belongs to $J$, but
\[
 (e_0f_0)^N\cdot1=(-e_0e_1)^N\ne0
 \quad\text{in }\gr_{\fm}Q\qquad(N\geq1).
\]
The direct-sum decomposition is compatible with the $\fm$-adic 
filtrations, so $J^N\gr_{\fm}M\ne0$ for every $N\geq1$.
Thus the weighted Casimir condition cannot be replaced by the same
coordinate ideal acting on the $\fm$-adicly graded module, even for
smooth mod-$p$ principal series.

\subsection{Filtered modules and Gelfand--Kirillov dimension}
\label{subsec:filtrations-dimension}

We write $\Krdim B$ for the Krull dimension of a commutative
noetherian ring $B$. For a finitely generated $B$-module $N$, put
\[
 \Krdim_B N:=\Krdim\Supp_{\Spec B}(N).
\]

Modules are left modules unless stated otherwise. We use decreasing,
exhaustive, separated filtrations indexed by $\Z$. Let $A$ have
ring filtration $F_A$. A compatible filtration $F_M$ on a finitely
generated $A$-module $M$ is good if there exist
$m_1,\ldots,m_t\in M$ and $d_1,\ldots,d_t\in\Z$ such that
\[
 F_M^nM=\sum_{i=1}^t(F_A^{n-d_i}A)m_i\qquad(n\in\Z).
\]
We write $\gr_{F_M}M$ for its associated graded module.
A filtered morphism $u:M\to M'$ is strict if
$u(F_M^nM)=u(M)\cap F_{M'}^nM'$ for every $n$. When considering
submodules and quotients of a filtered module, we use the induced
and quotient filtrations: for $N\subseteq M$,
\[
 F_N^nN=N\cap F_M^nM,\qquad
 F_{M/N}^n(M/N)=(F_M^nM+N)/N.
\]
After the reindexing $F_n=F^{-n}$ these are the conventions of
\cite[\S\S~I.4--I.5]{LvO96}.
Any two good filtrations differ by a bounded shift
\cite[Lemma~I.5.3]{LvO96}. The filtrations used below are complete.

Unless a different module filtration is specified, a finitely generated
$\Lambda$-module $M$ carries the natural filtration
\[
 F_\omega^nM=(F_\omega^n\Lambda)M,
 \qquad \gr_\omega M:=\gr_{F_\omega}M.
\]
This is the quotient filtration obtained from any finite generating
set placed in degree zero: each $F_\omega^n\Lambda$ is a two-sided
ideal, so the filtration is independent of that generating set.
It is complete, separated and good. We use the same convention for
finitely generated $\widetilde\Lambda$-modules, with the fixed
filtration $F_\omega^\bullet\widetilde\Lambda$.

\begin{lemma}\label{lem:good-filtration-annihilation}
Let $F$ and $F'$ be good filtrations on a finitely generated module
$M$ over the same filtered ring $A$, and let $\mathfrak b$ be a
homogeneous ideal of $\gr A$. If
$\mathfrak b^N\gr_F M=0$ for some $N\geq1$, then
$\mathfrak b^{N'}\gr_{F'}M=0$ for some $N'\geq1$.
In particular, whether some power of $J$ annihilates the associated
graded module is independent of the good module filtration, once
the ring filtration is fixed.
\end{lemma}

\begin{proof}
Choose $c\geq0$ with
$F^{n+c}M\subseteq F'^nM\subseteq F^{n-c}M$ for every $n$.
A lift $\widetilde b\in F_A^dA$ of a homogeneous
$b\in\mathfrak b^N$ of degree $d$ satisfies
$\widetilde b F^nM\subseteq F^{n+d+1}M$.
Take $\ell=2c+1$ such lifts, of degrees $d_1,\ldots,d_\ell$, and
put $D=d_1+\cdots+d_\ell$. Then
\[
 \widetilde b_1\cdots\widetilde b_\ell F'^nM
 \subseteq F^{n-c+D+\ell}M\subseteq F'^{n+D+1}M.
\]
These products span $\mathfrak b^{N\ell}$, so $N'=N(2c+1)$ works.
\end{proof}

For a nonzero finitely generated $\Lambda$-module $M$, define
\[
 j_\Lambda(M)=\min\{i\geq0:\Ext^i_\Lambda(M,\Lambda)\ne0\},
 \qquad \delta_H(M)=3n-j_\Lambda(M).
\]
This is the canonical Iwasawa dimension. It is unchanged on restriction
to an open uniform subgroup. For an admissible smooth $G$-representation
$\pi$ with a central character, put
$\GKdim(\pi)=\delta_H(\pi^\vee)$, with $\GKdim(0)=-\infty$.
We likewise put $\delta_H(0)=-\infty$.
This agrees with the canonical dimension of $\pi^\vee$ over
$\widetilde\Lambda$, since $Z_H$ acts trivially. See \cite[\S\S~5.4--5.5]{AB06},
\cite[Thm.~3.1 and Prop.~3.2]{Koh17}, and \cite[\S~6.1]{Sor26b}.

\begin{proposition}\label{prop:canonical-dimension}
For a nonzero finitely generated $\Lambda$-module $M$ and any good
filtration $F_M$ compatible with $F_\omega^\bullet\Lambda$,
\[
 \delta_H(M)=\Krdim_B\gr_{F_M}M.
\]
\end{proposition}

\begin{proof}
By Proposition~\ref{prop:graded-algebra}, $B$ is a polynomial ring
in $3n$ variables. The filtered comparison theorem
\cite[III.2.5, Thm.~2]{LvO96} and the grade--codimension formula
for polynomial rings give
\[
 j_\Lambda(M)=j_B(\gr_{F_M}M)
 =3n-\Krdim_B\gr_{F_M}M.
\]
The assertion follows from the definition of $\delta_H(M)$.
\end{proof}

\subsection{Microlocalization and control of ideals}\label{subsec:control}

We recall the filtered-algebra results used in the faithfulness proof,
using decreasing filtrations throughout. The conventions of
\cite{AW09} are translated by reindexing $F_n=F^{-n}$.

Let $A$ be a complete filtered $\F$-algebra such that
$B=\gr A$ is a commutative noetherian domain.
Write $\gr(a)$ for the principal symbol of a nonzero element
$a\in A$. For a multiplicative set $S$ of nonzero homogeneous
elements of $B$, write $Q_S(A)$ for the completion of the
Ore localization at
\[
 \widetilde S=\{a\in A\setminus\{0\}:\gr(a)\in S\}.
\]
Thus we first invert the elements whose principal symbols lie in
$S$, and then complete for the induced filtration.
For $0\ne b\in B$ homogeneous, we also write $Q_b(A)$ for
this microlocalization with $S=\{b^m:m\geq0\}$.
For a finitely generated left $A$-module $M$ with a good filtration, put
$Q_S(M)=Q_S(A)\otimes_{A}M$, with the tensor product filtration.
Then
\begin{equation}\label{eq:microlocalization}
 \gr Q_S(A)\cong S^{-1}B,\qquad
 \gr Q_S(M)\cong S^{-1}\gr M.
\end{equation}
Microlocalization is exact on finitely generated modules and preserves
the induced filtrations on submodules. This is
\cite[Prop.~4.3 and Lemma~4.4]{AWZ08}, applied to the opposite ring
for left modules. The required noetherianity of the Rees ring follows
from completeness and noetherianity of the associated graded ring;
see \cite[Prop.~II.2.2.1]{LvO96}.

\begin{definition}[{\cite[Def.~2.2]{AWZ08}}]\label{def:frobenius-pair}
A pair $(A_0,A_1)$ is a Frobenius pair if $A_0$ is a complete filtered
$\F$-algebra, $A_1$ is a closed subalgebra with the induced filtration,
$B_0=\gr A_0$ is a commutative noetherian domain, and
$B_1=\gr A_1\subseteq B_0$ satisfies
\[
 B_0^p\subseteq B_1,\qquad
 B_0=\bigoplus_{\alpha\in\{0,\ldots,p-1\}^t}B_1w^\alpha
\]
for some homogeneous $w_1,\ldots,w_t\in B_0$.
\end{definition}

For a Frobenius pair $(A_0,A_1)$, a \emph{source of derivations} is a
sequence $\boldsymbol\alpha=(\alpha_r)_{r\geq0}$ of
$\F$-algebra endomorphisms of $A_0$ with integers
$\theta_r,\theta'_r\geq1$ such that, for every $n$,
\begin{equation}\label{eq:source-shifts}
 (\alpha_r-1)(F^nA_0)\subseteq F^{n+\theta_r}A_0,\qquad
 (\alpha_r-1)(F^nA_1)\subseteq F^{n+\theta'_r}A_0,
 \qquad \theta'_r-\theta_r\longrightarrow\infty.
\end{equation}
For $a\in F^nA_0$, the class of $(\alpha_r-1)(a)$ in
$F^{n+\theta_r}A_0/F^{n+\theta_r+1}A_0$ depends only on the class
of $a$ in $F^nA_0/F^{n+1}A_0$. These maps define derivations
$D_{\boldsymbol\alpha,r}$ of $B_0$ of degree $\theta_r$.
For a family $\mathfrak S$ of sources and a graded ideal
$\mathfrak b\subseteq B_0$, put
\[
 \mathfrak b^{\mathfrak S}=\{b\in B_0:
 \forall\,\boldsymbol\alpha\in\mathfrak S,\quad
 D_{\boldsymbol\alpha,r}(b)\in\mathfrak b\text{ for }r\gg0\}.
\]
The threshold for $r$ may depend on $b$ and $\boldsymbol\alpha$.
A right ideal $I\subseteq A_0$ is $\mathfrak S$-invariant if
$\alpha_r(I)\subseteq I$ for all sufficiently large $r$, for each
source. These are the conventions of \cite[\S\S~3.1--3.4]{AW09}.

\begin{theorem}[{\cite[Thm.~3.5 and its proof]{AW09}}]
\label{thm:filtered-control}
Let $(A_0,A_1)$ be a Frobenius pair, let $\mathfrak S$ be a family of
sources of derivations, and let $I$ be an $\mathfrak S$-invariant right
ideal of $A_0$. Put $\mathfrak b=\gr I$. If
\[
 \Der_{B_1}(B_0)(\mathfrak b^{\mathfrak S})\subseteq\mathfrak b,
\]
then
\[
 I=(I\cap A_1)A_0,\qquad
 \gr(I\cap A_1)=\mathfrak b\cap B_1,\qquad
 \mathfrak b=(\mathfrak b\cap B_1)B_0.
\]
\end{theorem}

The first equality means that every element of $I$ is a finite sum
of products $xa$ with $x\in I\cap A_1$ and $a\in A_0$.
Here $I$ and $I\cap A_1$ carry the induced filtrations.
The last two identities follow from the proof of the cited theorem.
If $S_1=S\cap B_1$, then
$(Q_S(A_0),Q_{S_1}(A_1))$ is again a Frobenius pair, and sources
extend with the same shifts and localized derivations; see
\cite[Lemma~2.2 and Prop.~3.2]{AW09}.

\section{The ramified torus \texorpdfstring{$T$}{T}}
\label{sec:ramified-torus}

We construct a compact subgroup $T\subseteq H$ and prove finite
generation over its Iwasawa algebra. We then compare the canonical
Iwasawa dimension with the Krull dimension over $\Lambda_T=\F[[T]]$.

\subsection{Construction and strictness}

Recall that $\Pi=E_++\varpi_K E_-$. Put
\[
 \Theta=\begin{cases}
 \Pi,&\text{split case},\\
 \varpi_D,&\text{quaternionic case}.
 \end{cases}
\]
In either case $\Theta^2=\varpi_K$, so $K'=K(\Theta)$ is a ramified
quadratic extension of $K$. Define
\begin{equation}\label{eq:ramified-torus}
 T=\exp(p^2\cO_K\Theta),\qquad \Lambda_T=\F[[T]].
\end{equation}
Since the elements of $p^2\cO_K\Theta$ commute, the exponential
satisfies $\exp(X+Y)=\exp(X)\exp(Y)$ on this subgroup. It therefore
gives a group isomorphism from the additive group $p^2\cO_K\Theta$
to $T$, with inverse $\log$. Every $X\in p^2\cO_K\Theta$ has trace
zero in $K'/K$, so $\exp(X)$ has norm one. Thus $T$ is an open subgroup
of the principal norm-one units of $K'$, lies in $H$, and
$T\cong\Z_p^n$.

An ordered basis of $T$ is
\[
 g^T_{i,j}:=\exp(p^2\alpha^i\varpi_K^j\Theta)
 \qquad(0\leq i\leq f-1,\ 0\leq j\leq e-1).
\]
Consequently $\Lambda_T$ is a power series ring in the $n$ variables
$g^T_{i,j}-1$. Equip it with the weighted power series filtration
$F_\omega^\bullet \Lambda_T$
assigning weight $a_j$ to $g^T_{i,j}-1$, as in
\eqref{eq:weighted-ring-filtration}. We compare this filtration
with the one induced from $\Lambda$. Put
\[
 t_{k,j}=e_{k,j}+f_{k,j},\qquad
 B_T=\F[t_{k,j}:0\leq k\leq f-1,\ 0\leq j\leq e-1]\subset B.
\]

\begin{proposition}\label{prop:torus-strictness}
The associated graded algebra of $\Lambda_T$ is $B_T$, with
$\deg_\omega t_{k,j}=a_j$. In particular, $\Lambda_T\hookrightarrow\Lambda$ is
strict for the power series filtration assigning weight $a_j$
to $g^T_{i,j}-1$, and this filtration on $\Lambda_T$ is complete.
\end{proposition}

\begin{proof}
In the split case the exponential expansion gives
\[
 \gr_\omega(g^T_{i,j}-1)
   =\gr_\omega(u^+_{i,j}-1)+\gr_\omega(u^-_{i,j}-1).
\]
The higher-order terms in this expansion, including those in the
diagonal matrix entries, have weight greater than $a_j$. Multiplying by
$\sigma_k(\overline\beta_i)$ and summing over $i$ gives
$e_{k,j}+f_{k,j}$. In the quaternionic case,
$g^T_{i,j}=u^+_{i,j}$. Since
$\sigma_{k+f}(\overline\beta_i)=\sigma_k(\overline\beta_i)$ and
$\sigma_{k+f}(\overline\eta)=-\sigma_k(\overline\eta)$,
the definitions in \eqref{eq:quaternionic-graded-coordinates} give
\[
 t_{k,j}=e_{k,j}+f_{k,j}
   =\sum_{i=0}^{f-1}\sigma_k(\overline\beta_i)
                         \gr_\omega(u^+_{i,j}-1).
\]
Thus, in either case,
\[
 t_{k,j}=\sum_{i=0}^{f-1}\sigma_k(\overline\beta_i)
                    \gr_\omega(g^T_{i,j}-1).
\]
The matrix $(\sigma_k(\overline\beta_i))_{k,i}$ is invertible,
so the symbols $\gr_\omega(g^T_{i,j}-1)$ are also
algebraically independent.

Expand an element of $\Lambda_T$ as an ordered power series in the
$g^T_{i,j}-1$. Its lowest nonzero weighted part cannot vanish in $B$,
since the $t_{k,j}$ are algebraically independent. Its weighted order
therefore equals its induced order in $\Lambda$. This proves strictness
and $\gr_\omega \Lambda_T=B_T$. The filtration is complete and defines
the topology given by powers of the maximal ideal
$(g^T_{i,j}-1:i,j)$ in $\Lambda_T$.
\end{proof}

\subsection{The Casimir quotient}

\begin{proposition}\label{prop:general-casimir-freeness}
The natural map $B_T\longrightarrow R$ is injective, and
$R$ is a finite free $B_T$-module of rank $4^n$.
In particular, $\Krdim R=n$.
\end{proposition}

\begin{proof}
We first prove that
$R/(t_{k,j})_{k,j}R$ is finite-dimensional over $\F$.
Let $\mathfrak p$ be a prime containing $J$ and all $t_{k,j}$.
In the integral domain $B/\mathfrak p$,
we have $F_k(X)=-E_k(X)$, so
\[
 H_k(X)^2=4X E_k(X)^2.
\]
If $E_k\ne0$, then $H_k\ne0$ and
$2\operatorname{ord}_X(H_k)=1+2\operatorname{ord}_X(E_k)$,
which is impossible. Here $\operatorname{ord}_X$ is the least
exponent of $X$ with a nonzero coefficient. Thus $E_k=H_k=F_k=0$. Since every
$\mu_{k,j}$ is nonzero, all $e_{k,j},f_{k,j},h_{k,j}$ belong
to $\mathfrak p$.
The radical of $J+(t_{k,j})_{k,j}$ is therefore the
homogeneous maximal ideal of $B$. This proves the asserted
finite-dimensionality. Lift a homogeneous $\F$-basis of
$R/(t_{k,j})_{k,j}R$ to $R$. Graded Nakayama shows that these finitely
many lifts generate $R$ over $B_T$.

Setting all $e_{k,j}$ and $h_{k,j}$ equal to zero gives a surjection
$R\twoheadrightarrow\F[f_{k,j}]$ sending $t_{k,j}$ to $f_{k,j}$.
Its composite with $B_T\to R$ is an isomorphism onto this polynomial
ring. Hence $B_T\to R$ is injective and $\Krdim R=n$.

To prove freeness and compute the rank, use the auxiliary grading
\[
 \deg_{\mathrm{tot}} e_{k,j}=\deg_{\mathrm{tot}} f_{k,j}=\deg_{\mathrm{tot}} h_{k,j}=1.
\]
The Hilbert series below refer to this grading. Each of the $2n$
coefficients defining $J$ in \eqref{eq:general-casimir-ideal} is
homogeneous of degree $2$, and each $t_{k,j}=e_{k,j}+f_{k,j}$ has
degree $1$. Put
\[
 R_0=R/(t_{k,j})_{k,j}R\cong R\otimes_{B_T}\F.
\]
We have proved that $R_0$ is finite-dimensional over $\F$.
Thus the $2n$ Casimir coefficients, followed by the $n$ elements
$t_{k,j}$, form a homogeneous system of parameters in
$B$: there are $3n$ of them in a polynomial ring
with $3n$ variables, and their common quotient is finite-dimensional.
A polynomial ring is Cohen--Macaulay, so this system of parameters
is a regular sequence. This means that each element is a nonzero
divisor after quotienting by all preceding elements. In particular,
the coefficients defining $J$ form a regular sequence in
$B$, and the $t_{k,j}$ form a regular sequence on $R$.

Quotienting a graded ring by a homogeneous nonzero divisor of
degree $a$ multiplies its Hilbert series by $1-u^a$.
The $2n$ quadratic Casimir coefficients therefore give
\[
 \operatorname{Hilb}^{\mathrm{tot}}_{R}(u)
   =\frac{(1-u^2)^{2n}}{(1-u)^{3n}}
   =\frac{(1+u)^{2n}}{(1-u)^n}.
\]
Quotienting next by the $n$ linear elements $t_{k,j}$ gives
\[
 \operatorname{Hilb}^{\mathrm{tot}}_{R_0}(u)
   =(1-u)^n\operatorname{Hilb}^{\mathrm{tot}}_{R}(u)
   =(1+u)^{2n}.
\]

Choose a homogeneous $\F$-basis $\overline b_1,\ldots,\overline b_r$
of $R_0$, and lift it to homogeneous elements $b_1,\ldots,b_r$ of $R$.
Write $\delta_i=\deg_{\mathrm{tot}} b_i$. Graded Nakayama gives a surjective
graded $B_T$-linear map
\[
 \Phi:\bigoplus_{i=1}^r B_T(-\delta_i)\longrightarrow R
\]
sending the free generator of the $i$th summand to $b_i$;
the shift $B_T(-\delta_i)$ places that generator in degree $\delta_i$.
Since $\sum_i u^{\delta_i}=\operatorname{Hilb}^{\mathrm{tot}}_{R_0}(u)$, the
Hilbert series of the source is
\[
 \sum_{i=1}^r\frac{u^{\delta_i}}{(1-u)^n}
   =\frac{\operatorname{Hilb}^{\mathrm{tot}}_{R_0}(u)}{(1-u)^n}
   =\operatorname{Hilb}^{\mathrm{tot}}_{R}(u).
\]
In each degree, $\Phi$ is therefore a surjection between
finite-dimensional $\F$-vector spaces of equal dimension, hence
an isomorphism. It follows that $R$ is free over $B_T$, of rank
\[
 r=\dim_\F R_0=\operatorname{Hilb}^{\mathrm{tot}}_{R_0}(1)=2^{2n}=4^n.
\]
The auxiliary grading is used only in this commutative-algebra
argument; the filtration on $\Lambda$ is the one defined by $\omega$.
\end{proof}

\begin{corollary}\label{cor:general-casimir-finite-generation}
Let $M$ be a finitely generated $\Lambda$-module with a complete
separated good filtration $F_M$ compatible with
$F_\omega^\bullet\Lambda$. If
$J^N\gr_{F_M}M=0$ for some $N\geq1$, then $M$ is finitely
generated over $\Lambda_T$, and
\[
 \delta_H(M)=\Krdim_{\Lambda_T} M.
\]
\end{corollary}

\begin{proof}
Proposition~\ref{prop:general-casimir-freeness} and the successive
quotients $J^i/J^{i+1}$ show that $B/J^N$ is
finite over $B_T$. Hence $\gr_{F_M}M$ is finite over $B_T$.
By completeness, 
we can lift generators of $\gr_{F_M}M$ over $B_T$ to generators of
$M$ over $\Lambda_T$, and $F_M$ is a good $\Lambda_T$-module filtration.
Since $B/J^N$ is finite over $B_T$,
we have 
\[
 \delta_H(M)=\Krdim_{B}\gr_{F_M}M=\Krdim_{B_T}\gr_{F_M}M=\Krdim_{\Lambda_T} M
\]
by Proposition~\ref{prop:canonical-dimension}.
\end{proof}

\section{Faithfulness over the ramified torus}\label{sec:faithfulness}

In this section we prove faithfulness over $\Lambda_T$ for
infinite-dimensional finitely generated $\Lambda$-modules satisfying
the Casimir annihilation condition, using the ideal control theorem
recalled in Section~\ref{subsec:control}.

\subsection{The subgroup tower}\label{subsec:casimir-tower}

Put $\kappa=p/\varpi_K^e\in\cO_K^\times$.
For $s\geq0$ put
\[
 \mathfrak l_s=
 \begin{cases}
 p^2\varpi_K\cO_K H_{\mathrm{diag}}\oplus p^2\varpi_K^s\cO_K E_+
      \oplus p^2\varpi_K^{s+1}\cO_K E_-,&\text{split},\\
 p^2\varpi_K\cO_K\eta\oplus
      p^2\varpi_K^s\Theta\cO_L,&\text{quaternionic},
 \end{cases}
 \qquad H_s=\exp(\mathfrak l_s).
\]
Since $[\mathfrak l_s,\mathfrak l_s]\subseteq p\mathfrak l_s$,
the group $H_s$ is uniform. Its index in $H_{s-1}$ is
$p^{2f}$ for $s\geq1$. We use the subgroup
\[
 \Gamma=\begin{cases}
  \exp(p^2\varpi_K\cO_KH_{\mathrm{diag}}),&\text{split},\\
  \exp(p^2\varpi_K\cO_K\eta),&\text{quaternionic},
 \end{cases}
\]
which is contained in every $H_s$.

Note $H_0=H$. Give $\Lambda_0=\Lambda$ the filtration above
and $\Lambda_s=\F[[H_s]]$ its induced filtration.
For the tower and source calculations below we use the unrescaled
positive rational weights; multiplying all weights and shifts by
$\ell_\omega$ gives the integer filtration $F_\omega$. Thus
$\deg_\omega\gr_\omega(b)=\ell_\omega w(b)$ for the weight $w(b)$
defined below.
The positive weights on $e_{k,j},f_{k,j}$ and on $h_{k,j}$
in $B_0=\gr_\omega\Lambda_0=B$, for $0\leq j\leq e-1$, are
\[
 c+j/e,\qquad c+(j+1/2)/e,
\]
respectively, where $c=1+1/(2e)-C$.
By \eqref{eq:valuation-constant-range}, $c>1/(p-1)$. Define
\[
 \xi(v)=p^{\lfloor v/e\rfloor}
       \bigl(c+(v\bmod e)/e\bigr),\qquad
 \xi_{\mathrm{h}}(v)=p^{\lfloor v/e\rfloor}
       \bigl(c+(v\bmod e+1/2)/e\bigr)
       \quad(v\geq0).
\]
Write $w(b)$ for the unrescaled filtration weight of $b\in\Lambda$,
as distinct from the group $p$-valuation $\omega$.
For $v=me+j$ with $0\leq j\leq e-1$, the weights $\xi(v)$ and
$\xi_{\mathrm h}(v)$ apply to $e_{k,j}^{p^m},f_{k,j}^{p^m}$ and
to $h_{k,j}^{p^m}$, respectively. Equivalently, for $v\geq0$ and
$a\in\cO_K^\times$, in the split case we have
\[
 \begin{aligned}
 \xi(v)&=w\bigl(\exp(p^2\varpi_K^v aE_+)-1\bigr)
       =w\bigl(\exp(p^2\varpi_K^{v+1}aE_-)-1\bigr),\\
 \xi_{\mathrm h}(v)
 &=w\bigl(\exp(p^2\varpi_K^{v+1}aH_{\mathrm{diag}})-1\bigr).
 \end{aligned}
\]
In the quaternionic case, for $a\in\cO_K^\times$ and
$b\in\cO_L^\times$,
\[
 \begin{aligned}
 \xi(v)&=w\bigl(\exp(p^2\varpi_K^v\Theta b)-1\bigr),\\
 \xi_{\mathrm h}(v)
 &=w\bigl(\exp(p^2\varpi_K^{v+1}a\eta)-1\bigr).
 \end{aligned}
\]
Both $\xi$ and $\xi_{\mathrm h}$
are strictly increasing functions on $\Z_{\geq0}$, as is immediate
from these formulas.

Set
\[
 m_j(s)=\max\{0,\lceil(s-j)/e\rceil\}.
\]
\begin{proposition}\label{prop:casimir-tower-graded}
For every $s\geq0$, the induced filtration on $\Lambda_s$ is complete,
and its associated graded ring is the following subring
of $B_0=B$:
\begin{equation}\label{eq:casimir-tower-graded}
 B_s:=\gr_\omega\Lambda_s
 =\F[h_{k,j},e_{k,j}^{p^{m_j(s)}},f_{k,j}^{p^{m_j(s)}}:k,j].
\end{equation}
Here $0\leq k\leq f-1$ and $0\leq j\leq e-1$.
Each pair $(\Lambda_s,\Lambda_{s+1})$ is a Frobenius pair.
At the step from $B_s$ to $B_{s+1}$, the $2f$ generators involving
$e_{k,j},f_{k,j}$ with $j=s\bmod e$ are replaced by their $p$th
powers; all other generators remain unchanged.
In particular,
\begin{equation}\label{eq:casimir-full-cycle}
 B_{s+e}
 =\F[h_{k,j},e_{k,j}^{p^{m_j(s)+1}},f_{k,j}^{p^{m_j(s)+1}}:k,j].
\end{equation}
\end{proposition}

\begin{proof}
Since $K/K_0$ is totally ramified of degree $e$, we have
$\cO_K=\bigoplus_{j=0}^{e-1}\varpi_K^j\cO_{K_0}$.
For a nonzero element $\sum_j\lambda_j\varpi_K^j$, with
$\lambda_j\in\cO_{K_0}$, the nonzero summands have distinct
valuations modulo $\Z$. Hence
\[
 v_p\!\left(\sum_j\lambda_j\varpi_K^j\right)
   =\min_{\lambda_j\ne0}\left(v_p(\lambda_j)+\frac je\right).
\]
Membership in $\varpi_K^s\cO_K$ is therefore equivalent to
$\lambda_j\in p^{m_j(s)}\cO_{K_0}$ for every $j$. Thus
\[
 \varpi_K^s\cO_K
   =\bigoplus_{j=0}^{e-1}p^{m_j(s)}\varpi_K^j\cO_{K_0}.
\]
In the quaternionic case we also use
$\cO_L=\cO_K\oplus\eta\cO_K$. This gives the following exponentials
of a $\Z_p$-basis of $\mathfrak l_s$:
\[
 (u^+_{i,j})^{p^{m_j(s)}},\quad
 (u^-_{i,j})^{p^{m_j(s)}},\quad u^0_{i,j},
\]
where $0\leq i\leq f-1$ and $0\leq j\leq e-1$.
Since $H_s$ is uniform, these elements form an ordered basis of $H_s$.

For an original basis element $g$ of $H$ and any $m\geq0$,
characteristic $p$ gives $g^{p^m}-1=(g-1)^{p^m}$. Since $B_0$ is
a domain, taking principal symbols gives
\[
 \gr_\omega(g^{p^m}-1)=\gr_\omega(g-1)^{p^m}.
\]
The changes of coordinates in
\eqref{eq:split-graded-coordinates} and
\eqref{eq:quaternionic-graded-coordinates} are invertible.
Raising their coefficients to the $p^{m_j(s)}$th power preserves
invertibility.
Consequently, the symbols of the above basis elements minus $1$
generate exactly the polynomial subring on the right of
\eqref{eq:casimir-tower-graded}, and these symbols are algebraically
independent by Proposition~\ref{prop:graded-algebra}.
This proves \eqref{eq:casimir-tower-graded}.
The remaining assertions are immediate.
\end{proof}

To compute the derivations induced by conjugation, define the
following principal symbols in $B_s$, for $0\leq j\leq e-1$.
In the split case put
\[
 \begin{aligned}
 e_{k,j}^{[s]}&=\sum_{i=0}^{f-1}\sigma_k(\overline\beta_i)
    \gr_\omega\bigl(\exp(p^2\alpha^i\varpi_K^{s+j}E_+)-1\bigr),\\
 f_{k,j}^{[s]}&=\sum_{i=0}^{f-1}\sigma_k(\overline\beta_i)
    \gr_\omega\bigl(\exp(p^2\alpha^i\varpi_K^{s+j+1}E_-)-1\bigr).
 \end{aligned}
\]
In the quaternionic case, for $0\leq k\leq 2f-1$, put
\[
 x_{k,j}^{[s]}=\frac12\sum_{i=0}^{f-1}\left(
 \sigma_k(\overline\beta_i)
    \gr_\omega\bigl(\exp(p^2\alpha^i\varpi_K^{s+j}\Theta)-1\bigr)
 +\sigma_k(\overline\beta_i/\overline\eta)
    \gr_\omega\bigl(\exp(p^2\alpha^i\varpi_K^{s+j}\eta\Theta)-1\bigr)
 \right),
\]
and also set
\[
 e_{k,j}^{[s]}=x_{k,j}^{[s]},\qquad
 f_{k,j}^{[s]}=x_{k+f,j}^{[s]}
 \qquad(0\leq k\leq f-1).
\]
Directly from the definitions, in either case, for $s\geq0$ and
$0\leq j\leq e-2$,
\[
 e_{k,j+1}^{[s]}=e_{k,j}^{[s+1]},\qquad
 f_{k,j+1}^{[s]}=f_{k,j}^{[s+1]}.
\]
The symbols $e_{k,j}^{[s]},f_{k,j}^{[s]}$ have weight $\xi(s+j)$.
At $s=0$, we have $e_{k,j}^{[0]}=e_{k,j}$ and $f_{k,j}^{[0]}=f_{k,j}$;
also $x_{k,j}^{[0]}=x_{k,j}$ in the quaternionic case.

\begin{proposition}\label{prop:casimir-level-generators}
For every $s\geq0$, in either case, the symbols
$h_{k,j},e_{k,j}^{[s]},f_{k,j}^{[s]}$ form a set of
algebraically independent polynomial generators of $B_s$:
\begin{equation}\label{eq:casimir-level-generators}
 B_s=\F[h_{k,j},e_{k,j}^{[s]},f_{k,j}^{[s]}:k,j],
\end{equation}
where $0\leq k\leq f-1$ and $0\leq j\leq e-1$.
Moreover, $B_s$ is a free $B_{s+1}$-module of rank $p^{2f}$,
with basis the monomials having exponents in $\{0,\ldots,p-1\}$
in $e_{k,0}^{[s]},f_{k,0}^{[s]}$.
\end{proposition}

\begin{proof}
Recall that $\kappa=p/\varpi_K^e\in\cO_K^\times$.
Write $s+j=me+j'$ with $0\leq j'\leq e-1$. Then $m\geq0$ and
\[
 m_{j'}(s)=\max\{0,\lceil(s-j')/e\rceil\}
       =\lceil m-j/e\rceil=m,\qquad
 \varpi_K^{s+j}=p^m\kappa^{-m}\varpi_K^{j'}.
\]
In the split case, the definitions of the coordinates give
\begin{equation}\label{eq:casimir-coordinate-carry}
 e_{k,j}^{[s]}=\sigma_k(\overline\kappa^{-m})
                  e_{k-m,j'}^{p^m},\qquad
 f_{k,j}^{[s]}=\sigma_k(\overline\kappa^{-m})
                  f_{k-m,j'}^{p^m}.
\end{equation}
For the first identity, expand its right-hand side. Since
$\sigma_{k-m}(a)^{p^m}=\sigma_k(a)$ for $a\in k_K$ and
$\gr_\omega(g^{p^m}-1)=\gr_\omega(g-1)^{p^m}$, we obtain
\[
 \sigma_k(\overline\kappa^{-m})e_{k-m,j'}^{p^m}
 =\sum_{i=0}^{f-1}\sigma_k(\overline\kappa^{-m}\overline\beta_i)
   \gr_\omega\bigl(
       \exp(p^2\kappa^m\alpha^i\varpi_K^{s+j}E_+)-1\bigr).
\]
Write $\operatorname{Tr}$ for $\operatorname{Tr}_{k_K/\mathbf F_p}$ and put
\[
 S_{i'}=\gr_\omega\bigl(
       \exp(p^2\alpha^{i'}\varpi_K^{s+j}E_+)-1\bigr).
\]
Trace duality gives
\[
 \overline\kappa^m\overline\alpha^i
   =\sum_{i'=0}^{f-1}a_{i',i}\overline\alpha^{i'},
 \qquad
 a_{i',i}=\operatorname{Tr}
   (\overline\beta_{i'}\overline\kappa^m\overline\alpha^i)
   \in\mathbf F_p.
\]
The upper unipotent parameter is additive. On taking principal
symbols, product cross terms and changes of parameter divisible by
$\varpi_K$ have greater weight. Thus the residue-field expansion yields
\[
 \gr_\omega\bigl(
       \exp(p^2\kappa^m\alpha^i\varpi_K^{s+j}E_+)-1\bigr)
   =\sum_{i'=0}^{f-1}a_{i',i}S_{i'}.
\]
The reverse dual basis expansion
$b=\sum_i\operatorname{Tr}(b\overline\alpha^i)\overline\beta_i$,
applied to $b=\overline\beta_{i'}\overline\kappa^m$, gives
\[
 \sum_{i=0}^{f-1}\overline\beta_i a_{i',i}
   =\overline\beta_{i'}\overline\kappa^m.
\]
Since $\sigma_k$ fixes $a_{i',i}\in\mathbf F_p$, the coefficient
of $S_{i'}$ in the preceding weighted sum is
\[
 \begin{aligned}
 \sum_{i=0}^{f-1}
   \sigma_k(\overline\kappa^{-m}\overline\beta_i)a_{i',i}
 &=\sigma_k\left(\overline\kappa^{-m}
                  \sum_{i=0}^{f-1}\overline\beta_i a_{i',i}\right)\\
 &=\sigma_k(\overline\kappa^{-m}
             \overline\beta_{i'}\overline\kappa^m)
  =\sigma_k(\overline\beta_{i'}).
 \end{aligned}
\]
The sum is therefore $\sum_{i'}\sigma_k(\overline\beta_{i'})S_{i'}
=e_{k,j}^{[s]}$, proving the first identity.
The lower coordinates satisfy the same identity.
Using the $2f$ residue-field embeddings in the quaternionic case
gives, in the same way,
\[
 x_{k,j}^{[s]}=\sigma_k(\overline\kappa^{-m})x_{k-m,j'}^{p^m}.
\]

As $j$ runs through $0,\ldots,e-1$, the residue $j'=(s+j)\bmod e$
runs through the same range exactly once. For each such $j$, the
map $k\mapsto k-m$ is a permutation of the
embedding indices. All the scalar factors are nonzero.
Thus the symbols $e_{k,j}^{[s]},f_{k,j}^{[s]}$ in
\eqref{eq:casimir-level-generators} are precisely nonzero scalar multiples
of $e_{k,j}^{p^{m_j(s)}},f_{k,j}^{p^{m_j(s)}}$ in
\eqref{eq:casimir-tower-graded}, with permuted indices; the
$h_{k,j}$ are unchanged. This proves both generation and algebraic
independence in \eqref{eq:casimir-level-generators}.
The asserted $B_{s+1}$-module basis follows from the description
of each step in Proposition~\ref{prop:casimir-tower-graded}.
\end{proof}

\subsection{Sources of derivations}\label{subsec:casimir-sources}

We construct sequences of conjugation automorphisms of $\Lambda_s$
and compute the derivations they induce on $B_s$. These sequences
will be used in the ideal control argument of
Section~\ref{subsec:casimir-local-control}.
Fix lifts $c_0,\ldots,c_{f-1}\in\cO_{K_0}^{\times}$ of an
$\F_p$-basis of the residue field. For $0\leq j\leq e-1$,
$c\in\{c_0,\ldots,c_{f-1}\}$, and $r\geq3$, put
\[
 \gamma_{j,c,r}=
 \begin{cases}
 \exp(p^r\varpi_K^j cH_{\mathrm{diag}}),&\text{split},\\
 \exp(p^r\varpi_K^j c\eta),&\text{quaternionic}.
 \end{cases}
\]
Since $\gamma_{j,c,r}\in\Gamma\subseteq H_s$ for every $s$,
conjugation defines an $\F$-algebra automorphism
\[
 \alpha^{[s]}_{j,c,r}:\Lambda_s\longrightarrow\Lambda_s,
 \qquad z\longmapsto\gamma_{j,c,r}z\gamma_{j,c,r}^{-1}.
\]

\begin{proposition}\label{prop:casimir-conjugation-sources}
For every $s\geq0$, there exists $r_0\geq3$ such that, for each
$0\leq j\leq e-1$ and $c\in\{c_0,\ldots,c_{f-1}\}$, the sequence
\[
 \boldsymbol\alpha^{[s]}_{j,c}
   :=(\alpha^{[s]}_{j,c,r})_{r\geq r_0}
\]
is a source of derivations for $(\Lambda_s,\Lambda_{s+1})$.
With the rational weight convention of Section~\ref{subsec:casimir-tower},
its filtration shifts on $\Lambda_s$ and $\Lambda_{s+1}$ are,
respectively, $\theta_{s,j}(r)$ and $\theta_{s+1,j}(r)$, where
\begin{equation}\label{eq:casimir-source-shifts}
 \theta_{s,j}(r)=p^r\xi(s+j)-\xi(s).
\end{equation}
The integer shifts in \eqref{eq:source-shifts} are obtained by
multiplying these quantities by $\ell_\omega$. In particular,
\[
 \theta_{s+1,j}(r)-\theta_{s,j}(r)\longrightarrow\infty.
\]
\end{proposition}

\begin{proof}
Fix the source parameters $j,c$ and $r$.
Let $g$ be one of the group elements whose symbols
$\gr_\omega(g-1)$ occur in the definitions of
$e_{k,j'}^{[s]},f_{k,j'}^{[s]}$, with
$0\leq j'\leq e-1$, and write $\gamma=\gamma_{j,c,r}$.
We first show that $\gamma g\gamma^{-1}g^{-1}-1$ has leading
weight $p^r\xi(s+j'+j)$.
This computes the leading term of the automorphism difference,
since
\[
 \alpha^{[s]}_{j,c,r}(g-1)-(g-1)
   =(\gamma g\gamma^{-1}g^{-1}-1)g,
\]
and $g$ has principal symbol $1$.

In the split case, put $v=s+j'$ and $t=p^r\varpi_K^jc$.
We first consider an upper unipotent generator
\[
 g=\exp(X),\qquad X=p^2\varpi_K^v aE_+,
 \qquad a\in\cO_K^\times.
\]
Conjugation by $\gamma$ multiplies the upper-right entry
$p^2\varpi_K^v a$ of $g$ by $\exp(2t)$.
Since $r\geq3$, the exponential series gives
$\exp(2t)-1=2t(1+\varepsilon)$ with $\varepsilon\in p\cO_K$.
The upper unipotent subgroup is abelian, so
\[
 u:=\gamma g\gamma^{-1}g^{-1}
 =\exp\bigl((\exp(2t)-1)X\bigr)
 =\exp(Z)\exp(\varepsilon Z),
\]
with
\[
 \qquad Z=2tX=2p^{r+2}\varpi_K^{v+j}acE_+.
\]
Since $2ac$ is a unit, the weight formulas of
Section~\ref{subsec:casimir-tower} give
\[
 w\bigl(\exp(Z)-1\bigr)=p^r\xi(v+j).
\]
On the other hand, $\varepsilon\in p\cO_K$ and the monotonicity
of $\xi$ imply
\[
 w\bigl(\exp(\varepsilon Z)-1\bigr)
       \geq p^{r+1}\xi(v+j)>p^r\xi(v+j).
\]
Thus $\gr_\omega(u-1)=\gr_\omega(\exp(Z)-1)$,
which gives both the leading term and its claimed weight.
The lower unipotent case is analogous, with
$X=p^2\varpi_K^{v+1}aE_-$ and
$Z=-2tX$; the weight is again
$p^r\xi(v+j)$.

In the quaternionic case, put $v=s+j'$ and
\[
 \begin{aligned}
 X&=p^2\varpi_K^v\Theta b_0,\qquad
 Y=-2p^r\varpi_K^jc\eta,\\
 Z&=XY=-2p^{r+2}\varpi_K^{v+j}\Theta b_0c\eta,
 \end{aligned}
\]
where $b_0\in\cO_L^\times$ is a basis coefficient.
Since $\gamma X\gamma^{-1}=X\exp(Y)$, set
\[
 \Delta=\gamma X\gamma^{-1}-X=X(\exp(Y)-1),\qquad
 \mathcal I=p^{r+3}\varpi_K^{v+j}\mathfrak j,
\]
where $\mathfrak j=\varpi_D\cO_D$ is the multiplicatively closed
lattice of Section~\ref{subsec:groups}. The exponential series gives
$\exp(Y)-1-Y\in p^{2r}\varpi_K^{2j}\cO_L$, so
\[
 X\in p^2\mathfrak j,\qquad
 \Delta\in p^{r+2}\varpi_K^{v+j}\mathfrak j,\qquad
 \Delta-Z\in\mathcal I.
\]
For $m\geq2$, every term in $(X+\Delta)^m-X^m$ contains at
least one factor $\Delta$. Multiplicative closure therefore gives
\[
 (X+\Delta)^m-X^m
   \in p^{r+2m}\varpi_K^{v+j}\mathfrak j.
\]
Since $v_p(m!)\leq m-1$, all these terms divided by $m!$ belong
to $\mathcal I$. Summing the exponential series yields
\[
 \exp(X+\Delta)-\exp(X)\equiv\Delta\pmod{\mathcal I}.
\]
Now $\exp(-X)-1\in p^2\mathfrak j$, so, writing
$u=\gamma\exp(X)\gamma^{-1}\exp(-X)$, we obtain
\[
 u-1\equiv\Delta\equiv Z\pmod{\mathcal I}.
\]
The same estimate gives $\exp(Z)-1\equiv Z\pmod{\mathcal I}$,
and hence $u\exp(-Z)\in1+\mathcal I$.
Both $u$ and $\exp(Z)$ have reduced norm one. The exponential
and logarithm argument of Section~\ref{subsec:groups}, applied to
$\mathcal I$, therefore gives
\[
 u\exp(-Z)\in
 \exp\bigl(\mathcal I\cap\ker\Trd\bigr)
 =\exp\bigl(p^{r+1}\varpi_K^{v+j}\mathfrak l\bigr).
\]
The ordered-basis argument of
Proposition~\ref{prop:casimir-tower-graded}, applied to this lattice,
gives weights $p^{r+1}\xi(v+j+j'')$ for the generators obtained
from the $\cO_L\Theta$ component and
$p^{r+1}\xi_{\mathrm h}(v+j+j'')$ for those obtained from the
$K\eta$ component, for $0\leq j''\leq e-1$.
Consequently, $u\exp(-Z)-1$ has weight at least
$p^{r+1}\xi(v+j)$, strictly greater than the weight
$p^r\xi(v+j)$ of $\exp(Z)-1$.
Thus $u-1$ and $\exp(Z)-1$ have the same principal symbol.
This proves the claimed leading term and its weight.

The shift relative to $g-1$ is therefore
\[
 p^r\xi(s+j'+j)-\xi(s+j').
\]
For sufficiently large $r$, the unique minimum over $0\leq j'\leq e-1$
occurs at $j'=0$, because $\xi$ is strictly increasing.
Conjugation fixes the diagonal generators. The identity
$\alpha(xy)-xy=(\alpha(x)-x)\alpha(y)+x(\alpha(y)-y)$
extends this minimum shift from the generators to ordered monomials,
and then to all of $\Lambda_s$ by continuity.
Thus the shift on $\Lambda_s$ is $\theta_{s,j}(r)$, and the same
calculation at level $s+1$ gives $\theta_{s+1,j}(r)$ on
$\Lambda_{s+1}$. Moreover,
\[
 \begin{split}
 \theta_{s+1,j}(r)-\theta_{s,j}(r)
 ={}&p^r\bigl(\xi(s+j+1)-\xi(s+j)\bigr)\\
    &-\bigl(\xi(s+1)-\xi(s)\bigr)\longrightarrow\infty.
 \end{split}
\]
These estimates verify \eqref{eq:source-shifts} for the sequence
$\boldsymbol\alpha^{[s]}_{j,c}$ after choosing $r_0$ sufficiently large.
\end{proof}

\begin{proposition}\label{prop:casimir-source-derivations}
For $s\geq0$ and sufficiently large $r$, let
$D^{[s]}_{j,c,r}$ be the derivation of $B_s$ induced by
$\alpha^{[s]}_{j,c,r}-1$ at the filtration shift $\theta_{s,j}(r)$
of Proposition~\ref{prop:casimir-conjugation-sources}, as in
Section~\ref{subsec:control}.
These derivations vanish on every $h_{k,j'}$ and on all generators
except $e_{k,0}^{[s]},f_{k,0}^{[s]}$ in either case.
For fixed $j$ and sufficiently large $r$, their $\F$-linear span
as $c$ ranges through $\{c_0,\ldots,c_{f-1}\}$ contains
\begin{equation}\label{eq:casimir-split-sources}
 (e_{k-r,j}^{[s]})^{p^r}\frac{\partial}{\partial e_{k,0}^{[s]}}
 -(f_{k-r,j}^{[s]})^{p^r}\frac{\partial}{\partial f_{k,0}^{[s]}}
 \qquad(0\leq k\leq f-1)
\end{equation}
in the split case, and
\begin{equation}\label{eq:casimir-quaternionic-sources}
 (x_{k-r,j}^{[s]})^{p^r}\frac{\partial}{\partial x_{k,0}^{[s]}}
 -(x_{k+f-r,j}^{[s]})^{p^r}
                 \frac{\partial}{\partial x_{k+f,0}^{[s]}}
 \qquad(0\leq k\leq 2f-1)
\end{equation}
in the quaternionic case.
\end{proposition}

\begin{proof}
By the proof of Proposition~\ref{prop:casimir-conjugation-sources},
the smallest shift occurs at coordinate index $0$ for sufficiently
large $r$. Conjugation fixes every $h_{k,j'}$, and the shifts for
positive coordinate indices are strictly larger than the minimum.
Hence the induced derivations
vanish on all the other generators listed in the statement.

Fix $j,c,r$ and abbreviate $D=D^{[s]}_{j,c,r}$.
In the split case, write
\[
 c\alpha^i=\sum_{\ell=0}^{f-1}\gamma_{i,\ell}\alpha^\ell,
 \qquad
 \gamma_{i,\ell}
 =\operatorname{Tr}_{K_0/\Q_p}(\beta_\ell c\alpha^i)\in\Z_p.
\]
By the leading-term calculation of
Proposition~\ref{prop:casimir-conjugation-sources},
\[
 \begin{aligned}
 D\left(\gr_\omega\bigl(
       \exp(p^2\alpha^i\varpi_K^sE_+)-1\bigr)\right)
 &=\gr_\omega\bigl(
       \exp(2p^{r+2}c\alpha^i\varpi_K^{s+j}E_+)-1\bigr)\\
 &=2\sum_{\ell=0}^{f-1}\overline\gamma_{i,\ell}
       \left[\gr_\omega\bigl(
       \exp(p^2\alpha^\ell\varpi_K^{s+j}E_+)-1\bigr)\right]^{p^r}.
 \end{aligned}
\]
Trace duality also gives
\[
 \beta_\ell c=\sum_{i=0}^{f-1}\beta_i\gamma_{i,\ell},
 \qquad
 \sum_{i=0}^{f-1}\sigma_k(\overline\beta_i)
                  \overline\gamma_{i,\ell}
     =\sigma_k(\overline\beta_\ell\,\overline c).
\]
Substituting into the definition of $e_{k,0}^{[s]}$ and
interchanging the sums therefore yields
\[
 \begin{aligned}
 D(e_{k,0}^{[s]})
 &=2\sum_{\ell=0}^{f-1}
       \left(\sum_{i=0}^{f-1}\sigma_k(\overline\beta_i)
                    \overline\gamma_{i,\ell}\right)
       \left[\gr_\omega\bigl(
       \exp(p^2\alpha^\ell\varpi_K^{s+j}E_+)-1\bigr)\right]^{p^r}\\
 &=2\sigma_k(\overline c)\sum_{\ell=0}^{f-1}
       \sigma_k(\overline\beta_\ell)
       \left[\gr_\omega\bigl(
       \exp(p^2\alpha^\ell\varpi_K^{s+j}E_+)-1\bigr)\right]^{p^r}\\
 &=2\sigma_k(\overline c)(e_{k-r,j}^{[s]})^{p^r}.
 \end{aligned}
\]
The last equality follows directly from the definition of
$e_{k-r,j}^{[s]}$ and
$\sigma_{k-r}(\overline\beta_\ell)^{p^r}=\sigma_k(\overline\beta_\ell)$.
For the lower unipotent subgroup, replace the parameters
$p^2\alpha^i\varpi_K^sE_+$ by
$p^2\alpha^i\varpi_K^{s+1}E_-$; the leading coefficient is
$-2$ instead of $2$, so the same calculation gives
\[
 D(f_{k,0}^{[s]})
   =-2\sigma_k(\overline c)(f_{k-r,j}^{[s]})^{p^r}.
\]
Let $U_k$ denote the derivation in
\eqref{eq:casimir-split-sources}. Then we have 
\[
 D^{[s]}_{j,c,r}=2\sum_{k=0}^{f-1}\sigma_k(\overline c)U_k.
\]

In the quaternionic case, write the input as
$\exp(p^2\varpi_K^s a\Theta)$. The leading change of its logarithm under conjugation is
$2p^{r+2}\varpi_K^{s+j}c\eta a\Theta$.
The residue-field basis
$(\overline\alpha^i,\overline\alpha^i\overline\eta)_{0\leq i\leq f-1}$
has trace-dual basis
$(\overline\beta_i/2,\overline\beta_i/(2\overline\eta))_{0\leq i\leq f-1}$.
Expand multiplication by $\overline c\,\overline\eta$ in this
basis and apply the same trace-dual reconstruction identity.
The preceding calculation with these $2f$ basis elements gives
\[
 D(x_{k,0}^{[s]})
   =2\sigma_k(\overline c\,\overline\eta)
                   (x_{k-r,j}^{[s]})^{p^r}.
\]
Let $V_k$ denote the derivation in
\eqref{eq:casimir-quaternionic-sources}. Pairing the indices $k$
and $k+f$, and using
$\sigma_{k+f}(\overline c)=\sigma_k(\overline c)$ and
$\sigma_{k+f}(\overline\eta)=-\sigma_k(\overline\eta)$, gives
\[
 D^{[s]}_{j,c,r}
   =2\sum_{k=0}^{f-1}\sigma_k(\overline c)
                         \sigma_k(\overline\eta)V_k.
\]
Here $V_{k+f}=-V_k$, so the first $f$ indices suffice.

Finally, as $c$ runs through $c_0,\ldots,c_{f-1}$, the matrix
$(\sigma_k(\overline c_i))_{0\leq k\leq f-1,\,0\leq i\leq f-1}$ is
invertible, because $(\overline c_i)$ is a $\mathbf F_p$-basis of
$k_K$ and the embedding matrix has the inverse furnished by the
trace-dual basis. In the split case, inverting this matrix isolates
each $U_k$. In the quaternionic case the coefficient matrix differs
only by the nonzero row factors $\sigma_k(\overline\eta)$, so it
similarly isolates each $V_k$. Since $2\ne0$ in $\F$, this proves
the two assertions about the spans.

To express the coefficients in the original coordinates, write
$s+j=me+j'$ with $0\leq j'\leq e-1$.
Applying \eqref{eq:casimir-coordinate-carry} to the coefficients in
\eqref{eq:casimir-split-sources} shows that the coefficients of
$\partial/\partial e_{k,0}^{[s]}$ and
$\partial/\partial f_{k,0}^{[s]}$ are
$\sigma_k(\overline\kappa^{-m})$ times
$e_{k-r-m,j'}^{p^{r+m}}$ and $-f_{k-r-m,j'}^{p^{r+m}}$,
respectively. The same scalar multiplies the two coefficients
in \eqref{eq:casimir-quaternionic-sources}, since
$\overline\kappa\in k_K$ implies
$\sigma_{k+f}(\overline\kappa^{-m})
=\sigma_k(\overline\kappa^{-m})$.
\end{proof}

\subsection{Microlocal control}
\label{subsec:casimir-local-control}

We use two families of homogeneous polynomials. In the split case,
for $a=(k_1,j_1)$ and $b=(k_2,j_2)$ let
$\rho_{ab}\in\{0,\ldots,f-1\}$ be congruent to $k_2-k_1$
modulo $f$, and put
\begin{equation}\label{eq:casimir-split-minors}
 \Delta_{ab}
 =e_{k_1,j_1}^{p^{\rho_{ab}}}f_{k_2,j_2}
       -f_{k_1,j_1}^{p^{\rho_{ab}}}e_{k_2,j_2}.
\end{equation}
Omit identically zero members. In the quaternionic case use
\begin{equation}\label{eq:casimir-quaternionic-minors}
 \Delta_{k,j}=x_{k,j}^{q+1}-x_{k+f,j}^{q+1}.
\end{equation}
Since $\deg_\omega e_{k,j}=\deg_\omega f_{k,j}$ and
$\deg_\omega x_{k,j}=\deg_\omega x_{k+f,j}$, each $\Delta_{ab}$ and
$\Delta_{k,j}$ is homogeneous.

Fix a nonzero $\Delta=\Delta_{ab}$ in the split case, or
$\Delta=\Delta_{k,j}$ in the quaternionic case.
For the graded subrings $B_s\subset B_0=B$, set
\[
 S=\{\Delta^m:m\geq0\},\qquad S_s=S\cap B_s,\qquad
 \Lambda_{s,\Delta}=Q_{S_s}(\Lambda_s),\qquad
 B_{s,\Delta}=S_s^{-1}B_s.
\]
Here $S_s$ is a set of homogeneous symbols, as in
Section~\ref{subsec:control}. The power
$\Delta^{p^{\lceil s/e\rceil}}$ belongs to $S_s$;
inverting it inverts every element of $S_s$.
Also $S_s\cap B_{s+1}=S_{s+1}$.
Proposition~\ref{prop:casimir-tower-graded} and the microlocalization
results \cite[Lemma~2.2 and Prop.~3.2]{AW09} give compatible
closed inclusions $\Lambda_{s+1,\Delta}\subseteq\Lambda_{s,\Delta}$,
\[
 \gr_\omega\Lambda_{s,\Delta}=B_{s,\Delta}
   =B_s[(\Delta^{p^{\lceil s/e\rceil}})^{-1}],
\]
and Frobenius pairs $(\Lambda_{s,\Delta},\Lambda_{s+1,\Delta})$.
For each $s$, let $\mathfrak S_s$ consist of the extensions of the
sources $\boldsymbol\alpha^{[s]}_{j,c}$ in
Proposition~\ref{prop:casimir-conjugation-sources} to
$(\Lambda_{s,\Delta},\Lambda_{s+1,\Delta})$, for $0\leq j\leq e-1$ and
$c\in\{c_0,\ldots,c_{f-1}\}$. These extensions retain the shifts
in \eqref{eq:casimir-source-shifts}.

\begin{proposition}\label{prop:casimir-relative-control}
For every $s\geq0$ and every graded ideal
$\mathfrak b\subseteq B_{s,\Delta}$, we have 
\[
 \Der_{B_{s+1,\Delta}}(B_{s,\Delta})
       (\mathfrak b^{\mathfrak S_s})\subseteq\mathfrak b.
\]
Consequently, for a two-sided ideal $I_0\subseteq\Lambda_{0,\Delta}$
and $I_s=I_0\cap\Lambda_{s,\Delta}$,
\begin{equation}\label{eq:casimir-controlled-ideals}
 \gr I_s=(\gr I_{s+1})B_{s,\Delta},\qquad
 \gr I_{s+1}=\gr I_s\cap B_{s+1,\Delta}.
\end{equation}
\end{proposition}

\begin{proof}
The algebra $B_{s,\Delta}$ is a free $B_{s+1,\Delta}$-module with basis
\[
 \left\{\prod_{k=0}^{f-1}
       \bigl(e_{k,0}^{[s]}\bigr)^{a_k}
       \bigl(f_{k,0}^{[s]}\bigr)^{b_k}
       : 0\leq a_k,b_k<p\right\}.
\]
Consequently, $\Der_{B_{s+1,\Delta}}(B_{s,\Delta})$ is generated,
as a $B_{s,\Delta}$-module, by the partial derivatives with respect to
$e_{k,0}^{[s]}$ and $f_{k,0}^{[s]}$. We show that these partial derivatives are
$B_{s,\Delta}$-linear combinations of those derivations in 
Proposition~\ref{prop:casimir-source-derivations}.

In the split case write $\Delta=\Delta_{ab}$, with
$a=(k_1,j_1)$ and $b=(k_2,j_2)$ as above.
Fix $0\leq k\leq f-1$. To recover
$\partial/\partial e_{k,0}^{[s]}$ and
$\partial/\partial f_{k,0}^{[s]}$ from
\eqref{eq:casimir-split-sources}, choose
$j'_1,j'_2\in\{0,\ldots,e-1\}$ so that
$s+j'_1=m_1e+j_1$ and $s+j'_2=m_2e+j_2$.
Choose an arbitrarily large integer $r_0$ with
$r_0\equiv k-k_2\pmod f$, and use source indices
\[
 r_1=r_0+\rho_{ab}-m_1,\qquad r_2=r_0-m_2.
\]
These indices tend to infinity with $r_0$. Write $U_k(j,r)$ for
the derivation in \eqref{eq:casimir-split-sources}, displaying
its dependence on $j,r$. By
Proposition~\ref{prop:casimir-source-derivations}, each $U_k(j,r)$
is an $\F$-linear combination of the source derivations
$D^{[s]}_{j,c,r}$. The coordinate formula
\eqref{eq:casimir-coordinate-carry} gives
\[
 \begin{aligned}
 (e_{k-r_1,j'_1}^{[s]})^{p^{r_1}}
 &=\sigma_k(\overline\kappa^{-m_1})e_{k_1,j_1}^{p^{r_0+\rho_{ab}}},\\
 (e_{k-r_2,j'_2}^{[s]})^{p^{r_2}}
 &=\sigma_k(\overline\kappa^{-m_2})e_{k_2,j_2}^{p^{r_0}},
 \end{aligned}
\]
with the same identities for $f$. Thus the two normalized
$U_k$ have the following coefficient matrix:
\[
 \begin{pmatrix}
 \sigma_k(\overline\kappa^{m_1})U_k(j'_1,r_1)\\
 \sigma_k(\overline\kappa^{m_2})U_k(j'_2,r_2)
 \end{pmatrix}
 =
 \begin{pmatrix}
 e_{k_1,j_1}^{p^{r_0+\rho_{ab}}}&-f_{k_1,j_1}^{p^{r_0+\rho_{ab}}}\\
 e_{k_2,j_2}^{p^{r_0}}&-f_{k_2,j_2}^{p^{r_0}}
 \end{pmatrix}
 \begin{pmatrix}
 \partial/\partial e_{k,0}^{[s]}\\
 \partial/\partial f_{k,0}^{[s]}
 \end{pmatrix}.
\]
The matrix has determinant $-\Delta^{p^{r_0}}$, a unit in
$B_{s,\Delta}$. Inverting it expresses both partial derivatives
as $B_{s,\Delta}$-linear combinations of $U_k(j'_1,r_1)$ and
$U_k(j'_2,r_2)$, hence of source derivations with arbitrarily
large indices.

For $\Delta=\Delta_{k_1,j_1}$ in the quaternionic case, fix
$0\leq k\leq f-1$, choose $0\leq j'_1\leq e-1$ with
$s+j'_1=m_1e+j_1$, and take arbitrarily large
$r_0\equiv k-k_1\pmod{2f}$. The source indices
$r_1=r_0-m_1$ and $r_2=r_0+f-m_1$ give, after rescaling
each row by a nonzero scalar, the coefficient matrix for
$\partial/\partial x_{k,0}^{[s]}$ and
$\partial/\partial x_{k+f,0}^{[s]}$:
\[
 \begin{pmatrix}
 x_{k_1,j_1}^{p^{r_0}}&-x_{k_1+f,j_1}^{p^{r_0}}\\
 x_{k_1+f,j_1}^{p^{r_0+f}}&-x_{k_1,j_1}^{p^{r_0+f}}
 \end{pmatrix}.
\]
Its determinant is $-\Delta_{k_1,j_1}^{p^{r_0}}$, so the matrix is
invertible in $B_{s,\Delta}$ and recovers both
$\partial/\partial x_{k,0}^{[s]}$ and
$\partial/\partial x_{k+f,0}^{[s]}$. Letting $k$ vary gives every
generator of $\Der_{B_{s+1,\Delta}}(B_{s,\Delta})$ in either case.

For a graded ideal $\mathfrak b$ and an element
$z\in\mathfrak b^{\mathfrak S_s}$, choose the indices in these
two-row matrices large enough that each source derivation used
in the rows sends $z$ into $\mathfrak b$. The invertible matrices then show that
every relative partial derivative of $z$ belongs to $\mathfrak b$.
This is the derivation hypothesis of Theorem~\ref{thm:filtered-control}.
Each ideal $I_s=I_0\cap\Lambda_{s,\Delta}$ is preserved by the conjugation
sources: their conjugating elements belong to $\Gamma\subset H_s$
and are units of $\Lambda_{s,\Delta}$. Applying that theorem at each step proves
\eqref{eq:casimir-controlled-ideals}.
\end{proof}

The following corollary applies when an annihilator has a symbol that is
a polynomial in $e_{k,j},f_{k,j}$ with no $h_{k,j}$ variables.

\begin{corollary}\label{cor:casimir-ef-extraction}
Let $I_0$ be as in Proposition~\ref{prop:casimir-relative-control}.
If
\[
 \gr I_0\cap\F[e_{k,j},f_{k,j}:k,j]\ne0,
\]
then $I_0=\Lambda_{0,\Delta}$.
\end{corollary}

\begin{proof}
Denote the $2n$ variables $e_{k,j},f_{k,j}$ by
$X_1,\ldots,X_{2n}$, and write
$\boldsymbol h=(h_{k,j})_{k,j}$. At indices $s=em$, with $m\geq0$,
\[
 B_{em,\Delta}
 =\F[\boldsymbol h,X_1^{p^m},\ldots,X_{2n}^{p^m},(\Delta^{p^m})^{-1}].
\]
The ring $B_{0,\Delta}$ is free over $B_{em,\Delta}$ on
$X^\alpha$ with $0\leq\alpha_i\leq p^m-1$. This follows from the
polynomial-ring basis by localizing at $\Delta^{p^m}$.
Iterating \eqref{eq:casimir-controlled-ideals} gives
$\gr I_0=(\gr I_{em})B_{0,\Delta}$.
For a nonzero $\Psi(X)\in\gr I_0$, choose $p^m$ larger than
each individual variable degree of $\Psi$. In its unique expansion
in the basis $\{X^\alpha:0\leq\alpha_i\leq p^m-1\}$, each basis
coefficient is the coefficient of the same monomial in $\Psi$,
and therefore lies in $\F$. The equality
$\gr I_0=(\gr I_{em})B_{0,\Delta}$ forces every nonzero monomial
coefficient of $\Psi$ to belong to $\gr I_{em}$. Thus $1\in\gr I_0$. An element of $I_0$ with
principal symbol $1$ is a unit by completeness, proving the claim.
\end{proof}

We also need control of ideals in the Iwasawa algebras of the
upper and lower unipotent subgroups in the split case.
Let
\[
 N^+=\exp(p^2\cO_KE_+),\qquad
 N^-=\exp(p^2\varpi_K\cO_KE_-).
\]
Put $\Lambda^\pm:=\F[[N^\pm]]$, with the filtrations induced from $\Lambda$.
For $0\leq k_0\leq f-1$ and $0\leq j_0\leq e-1$, put $X=e_{k_0,j_0}$
and define $\Lambda_X^+:=Q_X(\Lambda^+)$, the microlocalization at
the powers of $X$. Its graded ring is
\[
 B_X^{+}:=\gr_\omega\Lambda_X^{+}
             =\F[e_{k,j}:k,j][X^{-1}].
\]
Similarly, for $X=f_{k_0,j_0}$, put $\Lambda_X^-:=Q_X(\Lambda^-)$.
Its graded ring is
\[
 B_X^{-}:=\gr_\omega\Lambda_X^{-}
             =\F[f_{k,j}:k,j][X^{-1}].
\]

For unramified $K$, an analogous result is proved in
\cite[Cor.~3.7]{BHH25}: the microlocalization of $\F[[(\cO_K,+)]]$
obtained by inverting all the principal symbols
$y_0,\ldots,y_{f-1}$ has no proper nonzero ideals stable under an
open subgroup of $\cO_K^\times$ acting by multiplication on $\cO_K$.
The lemma below allows arbitrary ramification and inverts only one
$e_{k_0,j_0}$ or $f_{k_0,j_0}$.

\begin{lemma}\label{lem:casimir-unipotent-rigidity}
For every $0\leq k_0\leq f-1$ and $0\leq j_0\leq e-1$, the only
$\Gamma$-stable ideals of $\Lambda_{e_{k_0,j_0}}^{+}$ or
$\Lambda_{f_{k_0,j_0}}^{-}$ are zero and the whole ring.
\end{lemma}

\begin{proof}
Fix $(k_0,j_0)$. For $s\geq0$, set
\[
 N_s^+=\exp(p^2\varpi_K^s\cO_KE_+),\qquad
 N_s^-=\exp(p^2\varpi_K^{s+1}\cO_KE_-),\qquad
 \Lambda_s^\pm:=\F[[N_s^\pm]].
\]
For $X=e_{k_0,j_0}$ in the upper case and $X=f_{k_0,j_0}$ in the
lower case, define
\[
 \Lambda_{s,X}^\pm:=Q_{X^{p^{m_{j_0}(s)}}}(\Lambda_s^\pm),
 \qquad m_j(s)=\max\{0,\lceil(s-j)/e\rceil\}.
\]
The microlocalization inverts the powers of the indicated element. 
We consider $\Lambda_X^+$ with $X=e_{k_0,j_0}$, 
and the case of $\Lambda_X^-$ with $X=f_{k_0,j_0}$ is the same. Its tower has graded rings
\[
 B_{s,X}^+:=\gr_\omega\Lambda_{s,X}^+
 =\F[e_{k,j}^{p^{m_j(s)}}:k,j]
       [(X^{p^{m_{j_0}(s)}})^{-1}].
\]
As in Proposition~\ref{prop:casimir-relative-control}, these
rings form Frobenius pairs $(\Lambda_{s,X}^+,\Lambda_{s+1,X}^+)$, and
$\Der_{B_{s+1,X}^+}(B_{s,X}^+)$ is generated by
$\partial/\partial e_{k,0}^{[s]}$.

Restrict the conjugation sources to $\Lambda_s^+$ and extend
them to $\Lambda_{s,X}^+$, denoting their family by $\mathfrak S_s$.
By Proposition~\ref{prop:casimir-source-derivations}, their
induced derivations span
$(e_{k-r,j}^{[s]})^{p^r}\partial/\partial e_{k,0}^{[s]}$.
For each $k$, choose $0\leq j\leq e-1$ with $s+j=m_1e+j_0$ and arbitrarily
large $r\equiv k-k_0-m_1\pmod f$. The coordinate formula
\eqref{eq:casimir-coordinate-carry} gives
\[
 (e_{k-r,j}^{[s]})^{p^r}
 =\sigma_k(\overline\kappa^{-m_1})X^{p^{r+m_1}},
\]
which is a unit in $B_{s,X}^+$ for sufficiently large $r$.
Thus every relative partial derivative is a $B_{s,X}^+$-linear
combination of source derivations with arbitrarily large indices.
The same argument as in
Proposition~\ref{prop:casimir-relative-control} verifies the
derivation hypothesis.

Let $I\subseteq\Lambda_X^+$ be a $\Gamma$-stable ideal and set
$I_s=I\cap\Lambda_{s,X}^+$. Since the $\Lambda_{s,X}^+$ are preserved by $\Gamma$,
each $I_s$ is invariant under $\mathfrak S_s$.
Theorem~\ref{thm:filtered-control} therefore gives
\[
 I_s=I_{s+1}\Lambda_{s,X}^+,\qquad
 \gr I_s=(\gr I_{s+1})B_{s,X}^+.
\]
In particular, $\gr I=(\gr I_{en})B_X^+$ for every $n\geq0$, with
\[
 B_{en,X}^+=\F[e_{k,j}^{p^n}:k,j][(X^{p^n})^{-1}].
\]
If $I\ne0$, clearing the denominator of a nonzero element of
$\gr I$ gives $0\ne\Psi\in\gr I\cap\F[e_{k,j}:k,j]$.
Choosing $p^n$ larger than each individual variable degree of
$\Psi$, the coefficient-extraction argument of
Corollary~\ref{cor:casimir-ef-extraction} gives
$1\in\gr I$, hence a unit in $I$ by completeness. Thus $I=\Lambda_X^+$, as desired. 
\end{proof}

\subsection{Faithfulness}\label{subsec:casimir-main}

We now use the preceding control results to prove faithfulness
over the ramified torus. The following lemma describes the locus
inside $\operatorname{V}(J)$ where all the polynomials $\Delta$ defined
in Section~\ref{subsec:casimir-local-control} vanish.

\begin{lemma}\label{lem:casimir-support-charts}
In $\Spec B$, we have 
\[
 \operatorname{V}(J)\cap\bigcap_{a,b}\operatorname{V}((\Delta_{ab}))
 =\operatorname{V}((h_{k,j},e_{k,j}:k,j))
       \ \cup\ \operatorname{V}((h_{k,j},f_{k,j}:k,j))
\]
in the split case, and 
\[
 \operatorname{V}(J)\cap\bigcap_{k,j}\operatorname{V}((\Delta_{k,j}))
 =\operatorname{V}((h_{k,j},e_{k,j},f_{k,j}:k,j))
 =\{\mathfrak m_0\}
\]
in the quaternionic case, 
where 
$\mathfrak m_0=(h_{k,j},e_{k,j},f_{k,j}:k,j)$ is the maximal ideal
of the origin.
\end{lemma}

\begin{proof}
Let $\mathfrak p$ be a prime ideal of $B$ containing
$J$ and all the polynomials $\Delta$.
Work in $\Frac(B/\mathfrak p)$.
In the split case, taking $a=(k,j)$ and $b=(k,j')$ gives
$e_{k,j}f_{k,j'}=f_{k,j}e_{k,j'}$ for every $j,j'$.
If $E_k(X)$ and $F_k(X)$ are both nonzero, these equalities
make them proportional: $E_k(X)=\lambda F_k(X)$ for some
nonzero scalar $\lambda$.
The Casimir equation would then give
$H_k(X)^2=-4\lambda X F_k(X)^2$, which is impossible by degree considerations. Thus for every $k$, either $E_k(X)=0$ or
$F_k(X)=0$, and the Casimir equation then gives $H_k(X)=0$.
If $e_{k,j}\ne0$ and $f_{k',j'}\ne0$ for distinct $k,k'$,
then $f_{k,j}=e_{k',j'}=0$, so $\Delta_{(k,j),(k',j')}\ne0$.
This proves the split assertion.

In the quaternionic case, vanishing of all $\Delta_{k,j}$ implies
$x_{k,j}=0$ if and only if $x_{k+f,j}=0$. The two polynomials
$E_k(X),F_k(X)$ therefore have the same
degree in $X$, which makes $E_k(X)F_k(X)=0$ by the same degree considerations.
Thus both polynomials vanish
for every $k$, and $H_k(X)=0$ as well. 

Conversely, the two indicated coordinate subspaces in the split
case and the origin in the quaternionic case satisfy the Casimir
equations and make all the relevant $\Delta$ vanish.
\end{proof}

We also need two filtered-algebra facts. 

\begin{lemma}\label{lem:microlocalized-annihilator}
Let $A$ be a complete filtered ring whose associated graded ring
is a commutative noetherian domain, and let $Q$ be a
microlocalization as in Section~\ref{subsec:control}.
If $M$ is a finite $A$-module with a good filtration, then
every element of $\Ann_{A}(M)$ annihilates
$Q\otimes_{A}M$.
\end{lemma}

\begin{proof}
Since $\gr A$ is commutative, the leading terms of $sz$ and $zs$
cancel. Hence
\[
 \deg([s,z])\geq\deg(s)+\deg(z)+1.
\]
For a denominator $s$, put $(\operatorname{ad}s)(z)=[s,z]$.
Repeated use of $as^{-1}=s^{-1}a+s^{-1}[s,a]s^{-1}$ gives the
finite identity
\[
 as^{-1}=\sum_{r=0}^{n}s^{-r-1}(\operatorname{ad}s)^r(a)
       +s^{-n-1}(\operatorname{ad}s)^{n+1}(a)s^{-1}.
\]
The $r$th summand has degree at least $\deg(a)-\deg(s)+r$;
the last term has degree at least $\deg(a)-\deg(s)+n+1$.
Both bounds tend to $+\infty$, so the series converges and the
remainder tends to zero. Letting $n\to\infty$ therefore gives
\begin{equation}\label{eq:annihilator-commutator-series}
 as^{-1}=\sum_{r\geq0}s^{-r-1}(\operatorname{ad}s)^r(a).
\end{equation}
Now take $a\in\Ann_{A}(M)$. If $z$ annihilates $M$, then
$[s,z]m=s(zm)-z(sm)=0$ for every $m\in M$.
Thus every $(\operatorname{ad}s)^r(a)$ also annihilates $M$, and
\[
 s^{-r-1}(\operatorname{ad}s)^r(a)\otimes m
 =s^{-r-1}\otimes(\operatorname{ad}s)^r(a)m=0.
\]
Applying the finite identity and letting its remainder tend to zero
shows that $a(s^{-1}\otimes m)=0$.
These fractions span the Ore-localized module, which is dense in
$Q\otimes_{A}M$. Multiplication by $a$ is continuous, so $a$
annihilates the whole microlocalized module.
\end{proof}

For localizations of the upper and lower unipotent subgroup
algebras in the split case, we also need to lift a graded
annihilator over a commutative ring.

\begin{lemma}\label{lem:commutative-filtered-annihilator}
Let $A$ be a commutative filtered ring such that
$\gr A$ is a noetherian domain. Let $M$ be a finite $A$-module
with a good filtration. If a nonzero homogeneous
element of $\gr A$ annihilates $\gr M$, then $\Ann_{A}(M)\ne0$.
\end{lemma}

\begin{proof}
Choose good filtered generators $m_1,\ldots,m_n$ of degrees
$e_1,\ldots,e_n$. Let $c$ be a nonzero homogeneous
annihilator of $\gr M$, of degree $d$, and choose a lift
$\widetilde c\in F^dA$. Since $c\gr M=0$, we can write
\[
 \widetilde c m_j=\sum_i a_{ij}m_i,\qquad
 a_{ij}\in F^{d+e_j-e_i+1}A.
\]
By the Cayley--Hamilton theorem,
$a:=\det(\widetilde c I_n-(a_{ij}))$ annihilates $M$.
It remains to prove that $a\ne0$.

The determinant has the term $\widetilde c^{\,n}$, whose principal
symbol is $c^n$. On expanding the determinant and then each factor,
every other term is, up to sign, of the form
\[
 \widetilde c^{\,n-r}\prod_{j\in S}a_{\pi(j),j},
\]
where $\varnothing\ne S\subseteq\{1,\ldots,n\}$, $r=|S|$, and
$\pi$ is a permutation of $S$. 
Since $\pi$ permutes $S$, we have
\[
 \sum_{j\in S}(e_j-e_{\pi(j)})=0.
\]
Consequently,
\[
 \begin{aligned}
 \deg\!\left(\widetilde c^{\,n-r}
             \prod_{j\in S}a_{\pi(j),j}\right)
 &\ge (n-r)d+\sum_{j\in S}(d+e_j-e_{\pi(j)}+1)\\
 &=nd+r\ge nd+1.
 \end{aligned}
\]
Thus every term other than $\widetilde c^{\,n}$ has strictly higher
degree, so the principal symbol of $a$ is $c^n$.
This is nonzero because $\gr A$ is a domain, and hence $a\ne0$.
\end{proof}

\begin{theorem}\label{thm:torus-faithfulness}
For either choice of $H$, let $M$ be an infinite-dimensional
finitely generated $\Lambda$-module. If
\[
 J^N\gr_\omega M=0
 \qquad\text{for some }N\geq1,
\]
then $\Ann_{\Lambda_T}(M)=0$.
\end{theorem}

\begin{proof}
The module $\gr_\omega M$ is a finite graded $B$-module with
support contained in $\operatorname{V}(J)$. Its support is not
equal to the origin, since otherwise $\gr_\omega M$, and hence $M$,
would be finite-dimensional. Every microlocalized module below carries
the filtration induced from $F_\omega^\bullet M$; its associated
graded module is denoted by $\gr_\omega$ as well.

Suppose that $0\ne b\in\Ann_{\Lambda_T}M$. Then we have 
\[
 \gr_\omega(b)=\Psi(e_{k,j}+f_{k,j}:k,j),
 \qquad 0\ne\Psi\in\F[t_{k,j}:k,j].
\]
Thus $\gr_\omega(b)$ lies in $\F[e_{k,j},f_{k,j}:k,j]\setminus\{0\}$.
For any $\Delta$ from \eqref{eq:casimir-split-minors} or
\eqref{eq:casimir-quaternionic-minors} such that
$(\gr_\omega M)[\Delta^{-1}]\ne0$, put
\[
 \Lambda_\Delta:=Q_\Delta(\Lambda)=\Lambda_{0,\Delta},\qquad
 M_\Delta=\Lambda_\Delta\otimes_{\Lambda}M.
\]
Then $M_\Delta$ is nonzero, with
$\gr_\omega M_\Delta=(\gr_\omega M)[\Delta^{-1}]$.
By Lemma~\ref{lem:microlocalized-annihilator}, $b$ annihilates
$M_\Delta$. The symbol of $b$ in $\gr_\omega\Lambda_\Delta$ is still
$\Psi(e_{k,j}+f_{k,j}:k,j)\ne0$, because $B\to B[\Delta^{-1}]$ is injective.
Apply Corollary~\ref{cor:casimir-ef-extraction} to
$I_0=\Ann_{\Lambda_\Delta}M_\Delta$. It gives
$I_0=\Lambda_\Delta$, a contradiction.
Consequently $(\gr_\omega M)[\Delta^{-1}]=0$ for every nonzero polynomial
$\Delta$ in \eqref{eq:casimir-split-minors} or
\eqref{eq:casimir-quaternionic-minors}.

In the quaternionic case, Lemma~\ref{lem:casimir-support-charts}
now forces the support to be the origin, already a contradiction.
This proves faithfulness in that case.

In the split case, the same lemma confines the support to the union of the
subspaces $\operatorname{V}(h_{k,j},f_{k,j}:k,j)$ and
$\operatorname{V}(h_{k,j},e_{k,j}:k,j)$.
Choose, for example, $(k_0,j_0)$ and $X=e_{k_0,j_0}$ with
$(\gr_\omega M)[X^{-1}]\ne0$, and set
\[
 \begin{aligned}
 \Lambda_X&=Q_X(\Lambda),\qquad
 B_X=\gr_\omega\Lambda_X=B[X^{-1}],\\
 M_X&=\Lambda_X\otimes_{\Lambda}M,\qquad
 \gr_\omega M_X=(\gr_\omega M)[X^{-1}].
 \end{aligned}
\]
The alternative $X=f_{k_0,j_0}$ is symmetric.
The algebra $\Lambda_X^{+}$ from Lemma~\ref{lem:casimir-unipotent-rigidity}
embeds strictly in $\Lambda_X$: this holds before microlocalization by
the ordered-basis calculation, and afterward by the injection
of the graded localization $B_X^{+}\hookrightarrow B_X$.
The support of $\gr_\omega M_X$ lies in $h_{k,j}=f_{k,j}=0$ for all
$k,j$. Therefore a power of
$(h_{k,j},f_{k,j}:k,j)$ annihilates $\gr_\omega M_X$, so
$\gr_\omega M_X$ is finite over $B_X^{+}$.
Lift a finite homogeneous generating set of $\gr_\omega M_X$ to $M_X$
and iterate in the complete filtration, as in
Corollary~\ref{cor:general-casimir-finite-generation}.
This makes $M_X$ a finitely generated $\Lambda_X^{+}$-module with a
good filtration.

The action
$\gamma(q\otimes m)=\operatorname{Ad}(\gamma)(q)\otimes\gamma m$
makes $\Ann_{\Lambda_X^{+}}M_X$ a $\Gamma$-stable ideal.
The polynomial $\Psi(e_{k,j}+f_{k,j}:k,j)$ annihilates
$\gr_\omega M_X$, and some power of $(f_{k,j}:k,j)$ also annihilates it.
The difference
$\Psi(e_{k,j}+f_{k,j}:k,j)-\Psi(e_{k,j}:k,j)$
belongs to $(f_{k,j}:k,j)$, so a power of
$\Psi(e_{k,j}:k,j)$ annihilates $\gr_\omega M_X$. 
By Lemma~\ref{lem:commutative-filtered-annihilator}, we have 
$\Ann_{\Lambda_X^{+}}M_X\ne0$. Lemma~\ref{lem:casimir-unipotent-rigidity} forces
this annihilator to be the whole ring, contradicting $M_X\ne0$.
\end{proof}

\section{The main theorem and its applications}\label{sec:conclusion}

Now we can prove our main theorem.

\begin{theorem}\label{thm:main-lambda}
Let $M$ be an infinite-dimensional finitely generated
$\Lambda$-module. Suppose that
\[
 J^N\gr_\omega M=0
 \qquad\text{for some }N\geq1.
\]
Then $M$ is a finitely generated faithful $\Lambda_T$-module, and
\[
 \delta_H(M)=n.
\]
\end{theorem}

\begin{proof}
Corollary~\ref{cor:general-casimir-finite-generation} gives finite
generation over $\Lambda_T$ and identifies $\delta_H(M)$ with
$\Krdim_{\Lambda_T} M$. Theorem~\ref{thm:torus-faithfulness} gives
$\Ann_{\Lambda_T}(M)=0$. Since $M$ is finite over $\Lambda_T$, its support is
$\operatorname{V}(\Ann_{\Lambda_T}(M))=\Spec \Lambda_T$. Therefore
\[
 \delta_H(M)=\Krdim_{\Lambda_T} M=\Krdim \Lambda_T=n,
\]
as desired. 
\end{proof}

\begin{corollary}\label{cor:main-gk-dimension}
Let $G=\GL_2(K)$ or $D^\times$, and let $\pi$ be an
infinite-dimensional admissible smooth $\F$-representation of
$G$ with a central character. Suppose that its dual, restricted
to $H$, satisfies the Casimir annihilation condition of
Theorem~\ref{thm:main-lambda}. Then $\pi^\vee$ is a finitely
generated faithful $\Lambda_T$-module, and $\pi|_T$ is admissible.
Moreover,
\[
 \GKdim(\pi)=n.
\]
\end{corollary}

\begin{proof}
The central character is trivial on $Z_H$. Since
$\widetilde H=H\times Z_H$ is compact open in $G$, admissibility makes
$M=\pi^\vee$ a finitely generated $\Lambda$-module.
Apply Theorem~\ref{thm:main-lambda}. Finite generation over $\Lambda_T$
is equivalent to admissibility of $\pi|_T$. Canonical dimension is unchanged by passage to
an open subgroup or by removing a central subgroup acting
trivially, so $\GKdim(\pi)=\delta_H(M)=n$.
\end{proof}

Finally we apply the theorem to Banach representations considered in \cite{Sor26b}. For a finite
extension $E/\Q_p$, let $\varpi_E$ be a uniformizer and put
$k_E=\cO_E/\varpi_E\cO_E$.

\begin{corollary}\label{cor:banach-gk-dimension}
Let $G=\GL_2(K)$ or $D^\times$, let $E/\Q_p$ be finite, and
let $\mathcal V$ be an infinite-dimensional admissible
$E$-Banach representation of $G$. Suppose that the locally
$\Q_p$-analytic vectors $\mathcal V^{\mathrm{la}}$ have an
infinitesimal character. Then
\[
 \GKdim(\mathcal V)=n.
\]
\end{corollary}

\begin{proof}
We may enlarge $E$ to contain the images of all embeddings
$L\hookrightarrow\overline{\Q}_p$, so that $k_L$ embeds into $k_E$;
this does not change Gelfand--Kirillov dimension
\cite[Lemma~6.1]{Sor26b}. We apply the preceding results with $\F=k_E$.

Let $\widetilde H=H\times Z_H$ be the compact open subgroup of
\eqref{eq:norm-one-realization}. Choose a $\widetilde H$-stable
open bounded lattice $\mathcal V^\circ$ and put
\[
 \overline M=(\mathcal V^\circ/\varpi_E\mathcal V^\circ)^\vee.
\]
Admissible Banach duality makes $\overline M$ a finitely generated
$\widetilde\Lambda$-module \cite[\S\S~2.1 and 2.3]{Sor26b}.
It is infinite-dimensional, since otherwise completeness would make
$\mathcal V^\circ$ finite over $\cO_E$ and $\mathcal V$ finite-dimensional.

Put
\[
 \widetilde B=\gr_\omega\widetilde\Lambda
   =B[z_{k,j}:0\leq k\leq f-1,\ 0\leq j\leq e-1].
\]
By \cite[proof of Thm.~6.2, (24)]{Sor26b} and the identification
of the Casimir ideal in Section~\ref{subsec:casimir-ideal},
there is a complete separated good $\widetilde\Lambda$-module
filtration $F$ on $\overline M$, compatible with
$F_\omega^\bullet\widetilde\Lambda$, such that
\[
 \mathfrak b^N\gr_F\overline M=0,
 \qquad \mathfrak b=J\,\widetilde B+(z_{k,j}:k,j),
\]
for some $N\geq1$.
Since $\widetilde B/\mathfrak b\cong R$ is finite over $B_T$,
so is $\widetilde B/\mathfrak b^N$, and hence $\gr_F\overline M$.
By completeness, as in the proof of
Corollary~\ref{cor:general-casimir-finite-generation},
$\overline M$ is finite over $\Lambda_T$ and $F$ is a good
$\Lambda_T$-module filtration, hence a good $\Lambda$-module filtration.
Since $J^N\gr_F\overline M=0$, Lemma~\ref{lem:good-filtration-annihilation}
gives $J^{N'}\gr_\omega\overline M=0$ for some $N'\geq1$.
Theorem~\ref{thm:main-lambda} therefore gives $\delta_H(\overline M)=n$.

Since $\gr_F\overline M$ is finite over $B_T$, its Krull dimensions
over $\widetilde B$ and $B_T$ agree. Proposition~\ref{prop:canonical-dimension}
and \cite[Prop.~6.1]{Sor26b} now give
\[
 \GKdim(\mathcal V)
   =\Krdim_{\widetilde B}\gr_F\overline M
   =\Krdim_{B_T}\gr_F\overline M
   =\delta_H(\overline M)=n.
\]
The dimension comparison uses the full algebra
$k_E[[\widetilde H]]$, so the argument applies without a central character.
\end{proof}

\newpage

\begin{cbunit}
\appendix

\makeatletter
\let\saved@seccntformat\@seccntformat
\def\@seccntformat#1{\appendixname~\csname the#1\endcsname.\quad}
\section[Canonical diagrams of  supersingular representations]
{Canonical diagrams of  supersingular \\ representations}
\label{app:canonical-diagrams}
\let\@seccntformat\saved@seccntformat
\makeatother

\begin{center}
by \textsc{Yongquan Hu}\footnote{Morningside Center of Mathematics,
Academy of Mathematics and Systems Science, Chinese Academy of Sciences,
University of the Chinese Academy of Sciences, Beijing, 100190, China.
\textit{E-mail:} \href{mailto:yhu@amss.ac.cn}{\texttt{yhu@amss.ac.cn}}.}
\end{center}

\subsection{Introduction}

Let $p>2$ and let $K/\Q_p$ be a finite extension of degree
$n=ef$, with ramification index $e$ and residue degree $f$.
In this appendix the coefficient field $\F$ is algebraically closed
of characteristic $p$. To an irreducible smooth $\F$-representation $\pi$
of $\GL_2(K)$ admitting a central character,
we have attached a canonical diagram
$(C_0(\pi),C_1(\pi),\mathrm{can})$ which determines $\pi$
\cite[Thm.~1.1]{Hu12}.\footnote{It is denoted by
$(D_0(\pi),D_1(\pi),\mathrm{can})$ in \cite{Hu12}.}\smallskip

It is a natural question to determine how large $C_1(\pi)$ can be.
If $K=\Q_p$, then $C_1(\pi)$ is simply $\pi^{I_1}$, where
$I_1$ is the upper pro-$p$ Iwahori subgroup of $\GL_2(K)$.
By \cite[Thm.~1.3]{Hu12} and \cite[Cor.~2.10]{Sch15}, $C_1(\pi)$ is finite-dimensional over $\F$ if and only if $\pi$ is finitely presented. This criterion
has been used by Schraen \cite{Sch15} for quadratic extensions
and Wu \cite{Wu21} in general to prove that smooth  
supersingular representations of $\GL_2(K)$ with a central character
are not of finite presentation when $K\ne\Q_p$.\footnote{If the local field has characteristic $p$,
this failure of finite presentation is established in \cite{Hu12} by
proving $C_1(\pi)=C_0(\pi)=\pi$.} As a consequence, $C_1(\pi)$ is infinite-dimensional over $\F$. 
For some time, the author expected $C_1(\pi)$ to have a smaller Gelfand--Kirillov dimension than that of $\pi$, until Benjamin Schraen
pointed out that this fails when $[K:\Q_p]=2$: a result of Timmins \cite{Tim26} implies that $C_1(\pi)$ cannot have Gelfand--Kirillov dimension $1$. \smallskip

We use the small subgroup $H$ and the Casimir ideal
$J\subseteq\gr_\omega\F[[H]]$ of the split case in the main body
of the paper. The notion of Gelfand--Kirillov dimension is also
the same as there. The main result of this appendix is as follows.

\begin{theorem}\label{app:thm-canonical-diagram}
Suppose that $K\ne\Q_p$. Let $\pi$ be an admissible  supersingular  
$\F$-representation of $\GL_2(K)$ with a central character.
{Suppose that $\pi^\vee|_H$ has a complete separated good filtration,
compatible with the specified filtration on $\F[[H]]$}, such that
\[
 J^a\gr_\omega(\pi^\vee)=0
 \qquad\text{for some }a\geq1.
\]
Then
\[
 \GKdim(\pi)=n,\qquad C_1(\pi)=C_0(\pi)=\pi.
\]
\end{theorem}

The Casimir annihilation condition, introduced in \cite{Sor26b}, already
gives $\GKdim(\pi)=n$ by Corollary~\ref{cor:main-gk-dimension}. 
When $K/\Q_p$ is unramified and $p\geq5$, its
relation with the usual graded condition on the pro-$p$ Iwahori subgroup 
in \cite[\S~5.3]{BHH23} is explained in
Section~\ref{subsec:casimir-ideal}. We do not replace the
specified weighted filtration by the augmentation filtration in
the ramified case.\smallskip

All duals below use the $\F$-linear duality recalled in
Section~\ref{subsec:groups}. Admissible
smooth representations of a compact $p$-adic Lie group correspond
to finitely generated modules over its completed group algebra;
this holds for the present coefficient field by
\cite[Prop.~1.9]{Koh17}. In particular,
their subquotients are admissible.
For every compact group occurring below, canonical dimension is
computed after restriction to an open uniform subgroup; we write
$\delta(N)$ for this dimension. It is
independent of this choice by
\cite[Thm.~3.1 and Prop.~3.2]{Koh17}; no torsion-freeness
assumption on $I_1/Z_1$ is needed.\smallskip

\paragraph*{Notation.} Let $\cO_K$ be the ring of integers
of $K$, let $\varpi_K$ be its uniformizer, and let $q=p^f$.
Put
\[
 G=\GL_2(K),\qquad \mathcal K=\GL_2(\cO_K),\qquad
 \mathcal K_r=1+\varpi_K^rM_2(\cO_K)\quad(r\geq1).
\]
Write $Z$ for the center of $G$, $I$ for the upper triangular
Iwahori subgroup, $I_1$ for its pro-$p$ radical, and
$Z_1=Z\cap I_1=Z\cap\mathcal K_1$. We use the matrices
\[
 s=\begin{pmatrix}0&1\\1&0\end{pmatrix},\qquad
 \Pi=\begin{pmatrix}0&1\\\varpi_K&0\end{pmatrix},\qquad
 t=s\Pi=\begin{pmatrix}\varpi_K&0\\0&1\end{pmatrix}.
\]

Let $\pi$ be a smooth admissible irreducible
$\F$-representation of $G$ with a central character. Choose an
irreducible $\mathcal KZ$-subrepresentation $\sigma\subseteq\pi$,
and set
\[
 \begin{split}
 P^+&=\left\{\begin{pmatrix}a&b\\0&1\end{pmatrix}:
       a\in\cO_K\setminus\{0\},\ b\in\cO_K\right\},\\
 I^+&=\langle P^+\cdot\sigma\rangle,\qquad I^-=\Pi I^+,\\
 C_1&=I^+\cap I^-,\qquad C_0=\langle\mathcal K\cdot C_1\rangle.
 \end{split}
\]
The spaces $I^\pm$, and hence $C_0$ and $C_1$, are independent
of $\sigma$ by \cite[Cor.~3.17]{Hu12}. The spaces $I^\pm$ are $I$-stable, $C_1$ is stable
under $I$ and $\Pi$, and $\pi=I^++I^-$. The central character
is trivial on $Z_1$: its image there is a finite $p$-group,
whereas $\F^\times$ has no nontrivial $p$-power torsion.

\paragraph*{Acknowledgements.}
Y.~H. thanks Benjamin Schraen for pointing out that $C_1(\pi)$ cannot have Gelfand--Kirillov dimension $1$ in the quadratic case, and for helpful discussions on the subject.  
Y.~H. is partially supported by the CAS Project for Young Scientists in
Basic Research (Grant No.~YSBR-033), the National Natural Science Foundation
of China (Grants Nos.~12288201 and~12425103), the National Center for
Mathematics and Interdisciplinary Sciences, and the Hua Loo-Keng Key
Laboratory of Mathematics, Chinese Academy of Sciences. 

\paragraph*{AI Disclosure.}
The original version of the main theorem was established by the author for unramified $K$. I then used GPT-6 Astra to extend the result to arbitrary $K$.

\subsection{An exact sequence from \texorpdfstring{\cite{Hu12}}{Hu}}

Put
\[
 \Lambda=\F[[I_1/Z_1]],\qquad
 \Lambda_1=\F[[\mathcal K_1/Z_1]],\qquad
 M=(I^+/C_1)^\vee,\qquad Q=(C_0/C_1)^\vee.
\]
For a $\Lambda$-module $N$, let $N^t$ denote the
$\Lambda_1$-module obtained by the action
\[
 h\cdot_t v=(t^{-1}ht)v
 \qquad(h\in\mathcal K_1/Z_1,\ v\in N),
\]
and define
\[
 \mathcal F(N)=\Lambda\otimes_{\Lambda_1}N^t.
\]
Here $t^{-1}\mathcal K_1t\subset I_1$ is open, so $N^t$ is
finitely generated if $N$ is. The ordinary and completed tensor
products coincide, since $\Lambda$ is finite free as a right
$\Lambda_1$-module.

\begin{lemma}\label{app:lem-exact}
The modules $M$ and $Q$ are finitely generated, and there is an
exact sequence of $\Lambda$-modules
\begin{equation}\label{app:eq-exact}
 0\longrightarrow\mathcal F(M)\longrightarrow M
   \longrightarrow Q\longrightarrow0.
\end{equation}
\end{lemma}

\begin{proof}
Finite generation follows from admissibility: the spaces in this
sequence are subquotients of $\pi|_{I_1}$, and $\Lambda$ is
noetherian.

By \cite[Lem.~3.4]{Hu12}, $C_0\subseteq I^+$ and $C_0\cap I^-=C_1$.
It follows that
\begin{equation}\label{app:eq-decomposition}
 (\pi/C_0)|_I=(I^+/C_0)\oplus(I^-/C_1).
\end{equation}
By \cite[\S~3.1.1, preceding Lem.~3.1]{Hu12}, the evaluation map
$\operatorname{Ind}_I^{\mathcal K}I^-\twoheadrightarrow\pi$
has kernel contained in $\operatorname{Ind}_I^{\mathcal K}C_1$.
The latter subspace has image $C_0$. Its image and the kernel
containment show that it is the entire inverse image of $C_0$.
Hence we get an isomorphism
\begin{equation}\label{app:eq-induced-isomorphism}
 \operatorname{Ind}_I^{\mathcal K}(I^-/C_1)
       \xrightarrow{\ \sim\ }\pi/C_0.
\end{equation}

Mackey's formula and the double coset decomposition
$\mathcal K=I\sqcup IsI$ give
\[
 \bigl(\operatorname{Ind}_I^{\mathcal K}(I^-/C_1)\bigr)|_I
 =I^-/C_1\oplus
       \operatorname{Ind}_{I\cap sIs^{-1}}^I(I^-/C_1)^s.
\]
The two summands are sent via
\eqref{app:eq-induced-isomorphism} respectively onto $I^-/C_1$
and $I^+/C_0$ in \eqref{app:eq-decomposition}. Indeed, the first
map is the natural inclusion, and the second has image in
$I^+/C_0$ because $sI^-=tI^+\subseteq I^+$. Consequently
\[
 I^+/C_0\cong\operatorname{Ind}_{I\cap sIs^{-1}}^I(I^-/C_1)^s.
\]
Let $\mathcal  H$ be the subgroup of diagonal matrices with
Teichm\"uller entries in $\cO_K^\times$. Then
$I\cap sIs^{-1}=\mathcal  H\mathcal K_1$ and $I=\mathcal  H I_1$.
Restricting to $I_1$ therefore gives
\[
 (I^+/C_0)|_{I_1}
 \cong\operatorname{Ind}_{\mathcal K_1}^{I_1}(I^-/C_1)^s
 \cong\operatorname{Ind}_{\mathcal K_1}^{I_1}(I^+/C_1)^t.
\]
The last identification is induced by $\Pi$, with twist action
$h\cdot_a v=(a^{-1}ha)v$. Dualizing
\[
 0\longrightarrow C_0/C_1\longrightarrow I^+/C_1
   \longrightarrow I^+/C_0\longrightarrow0
\]
proves \eqref{app:eq-exact}.
\end{proof}

\subsection{A purity result}
 
\begin{proposition}\label{app:prop-purity}
Let $H$ be a compact $p$-adic Lie group. Let $V$ be an admissible smooth $\F$-representation of 
$H$, with $n=\GKdim(V)$. There is a largest subrepresentation
$P_{n,H}(V)\subset V$ whose dual is pure of dimension $n$. Moreover, the following statements hold.
\begin{enumerate}
\item[(i)] If $V\neq0$, then $P_{n,H}(V)\neq0$. 
\item[(ii)] If $H'$ is an open subgroup of $H$, then $P_{n,H}(V)=P_{n,H'}(V)$.
\item[(iii)] If $V_1\subset V$ is a subrepresentation and
$\GKdim(V/V_1)<n=\GKdim(V)$, then $P_{n,H}(V_1)=P_{n,H}(V)$  as subspaces of  $V$.
\end{enumerate}
\end{proposition}
\begin{proof}
Put $N:=V^\vee$. By assumption $\delta(N)=n$. Let $T_{n-1}(N)$ be the sum of all submodules of dimension
strictly less than $n$ (with the convention $\delta(0)=-\infty$). It follows from \cite[Prop.~3.5(vi)]{Ven02} that $P_{n,H}(V):=(N/T_{n-1}(N))^{\vee}$ satisfies the required property. Assertion (i) follows from the definition of $\delta(\cdot)$; see \cite[Prop.~3.5(ii)]{Ven02}. 

For (ii), note that the formation of $T_{n-1}(\cdot)$ on finitely generated $\F[[H]]$-modules commutes with restriction to open subgroups; see \cite[Prop.~2.7(ii)]{Ven02}. 
 
(iii) Dualize the exact sequence $0\rightarrow  V_1\rightarrow V\rightarrow V/V_1\rightarrow 0$. 
By assumption, the kernel of $V^\vee\twoheadrightarrow V_1^\vee$ has dimension
strictly less than $n$, and so does the inverse image of $T_{n-1}(V_1^\vee)$  by exactness. 
Thus this inverse image is contained in $T_{n-1}(V^{\vee})$. Conversely $T_{n-1}(V^{\vee})$
  maps into $T_{n-1}(V_1^\vee)$. Therefore that inverse image is  exactly $T_{n-1}(V^\vee)$.
This proves (iii).
\end{proof} \smallskip

\begin{corollary}\label{cor:pure}
Let $\pi$ be an admissible smooth irreducible $\F$-representation of $G$ with a central character. Then for any open subgroup $H$ of $G$, $\pi^{\vee}$ is a pure $\F[[H]]$-module.
\end{corollary}
\begin{proof}
Let $n=\GKdim(\pi)$ and $W=P_{n,\mathcal K}(\pi)\subset \pi$ be the subrepresentation constructed in Proposition~\ref{app:prop-purity} with respect to the open subgroup $\mathcal K$. In particular, $W$ is nonzero and stable under $\mathcal K$. By Proposition~\ref{app:prop-purity}(ii), $W$ is also equal to $P_{n,I}(\pi)$. Since $\Pi$ normalizes $I$, $W$ is also stable under $\Pi$. Since $G$ is generated by $\mathcal K$ and $\Pi$, $W$ is a $G$-subrepresentation of $\pi$, hence $W=\pi$ by irreducibility of $\pi$. 

The above discussion shows that $\pi^{\vee}$ is a pure $\F[[\mathcal{K}]]$-module. By Proposition~\ref{app:prop-purity}(ii), this also holds for any open subgroup $H$ of $G$.
\end{proof}

\subsection{Multiplicity comparison}
\label{app:subsec-multiplicity}

Let $U$ be a uniform pro-$p$ group. Its augmentation ideal
$\fm_U$ defines a complete separated filtration on $\F[[U]]$, with
\[
 \gr\F[[U]]\cong\F[X_1,\ldots,X_s],\qquad s=\dim U.
\]
If $N\ne0$ is a finitely generated $\F[[U]]$-module of canonical
dimension $n>0$, its Hilbert--Samuel function has the form
\begin{equation}\label{app:eq-hilbert-samuel}
 \dim_{\F}N/\fm_U^\ell N
   =\frac{\mu_U(N)}{n!}\ell^n+O(\ell^{n-1})
   \qquad(\ell\longrightarrow\infty),
\end{equation}
where $\mu_U(N)>0$. In dimension zero, $N$ is finite-dimensional
and we put $\mu_U(N)=\dim_{\F}N$. For modules of dimension at
most $n$, write
\[
 \mu_{U,n}(N)=
 \begin{cases}
 \mu_U(N),&\delta(N)=n,\\
 0,&\delta(N)<n.
 \end{cases}
\]
In particular $\mu_{U,n}(0)=0$. It is well known that
$\mu_{U,n}$ is additive on short exact sequences of modules
of dimension at most $n$.

\begin{lemma}\label{app:lem-multiplicity-scale}
Let $N\ne0$ be a finitely generated $\F[[U]]$-module of
dimension $n$. Then
\begin{equation}\label{app:eq-multiplicity-scale}
 \mu_{U^p}(N)=p^n\mu_U(N).
\end{equation}
If $U'\subseteq U$ is any open uniform subgroup, then
\begin{equation}\label{app:eq-multiplicity-monotonicity}
 \mu_{U'}(N)\geq\mu_U(N).
\end{equation}
\end{lemma}

\begin{proof}
The second assertion follows from
$\fm_{U'}^\ell N\subseteq\fm_U^\ell N$ and
\eqref{app:eq-hilbert-samuel}. For the first, choose an ordered
uniform basis $g_1,\ldots,g_s$ of $U$ and put $x_i=g_i-1$.
Then the $g_i^p$ form a basis of $U^p$, and $g_i^p-1=x_i^p$.
Set $\mathfrak b=\fm_{U^p}\F[[U]]$.
As $U^p$ is normal in $U$, this is a two-sided ideal and
$\mathfrak b^\ell=\fm_{U^p}^{\ell}\F[[U]]$ for $\ell\geq0$.
Moreover,
\begin{equation}\label{app:eq-augmentation-comparison}
 \fm_U^{p\ell+s(p-1)}\subseteq\mathfrak b^\ell
       \subseteq\fm_U^{p\ell}.
\end{equation}
The second inclusion follows from $x_i^p\in\fm_U^p$.
To see the first inclusion, the initial ideal of $\mathfrak b^\ell$
contains $(X_1^p,\ldots,X_s^p)^\ell$. Every monomial of degree
at least $p\ell+s(p-1)$ belongs to this ideal. Lifting its
homogeneous terms successively proves the inclusion by
completeness and closedness of $\mathfrak b^\ell$.

Applying \eqref{app:eq-augmentation-comparison} to $N$ gives
\[
 \dim_{\F}N/\fm_U^{p\ell}N
 \leq\dim_{\F}N/\fm_{U^p}^{\ell}N
 \leq\dim_{\F}N/\fm_U^{p\ell+s(p-1)}N,
\]
proving \eqref{app:eq-multiplicity-scale}. If $n=0$, all three
terms are eventually $\dim_{\F}N$, so the same formula holds.
\end{proof}

For $r\geq1$, let
\[
 U_r=\mathcal K_rZ_1/Z_1\subseteq I_1/Z_1.
\]
For $r\geq e+1$, these groups are uniform and satisfy
\begin{equation}\label{app:eq-congruence-powers}
 U_r^p=U_{r+e}.
\end{equation}
Indeed, since $p>2$, we may identify $U_r$ with the level-$r$
congruence subgroup of $\SL_2(\cO_K)$. For $r\geq e+1$,
the logarithm identifies $U_r$ with the powerful $\Z_p$-Lie lattice
$\varpi_K^r\mathfrak{sl}_2(\cO_K)$. The $p$th-power map
corresponds to multiplication by $p$, and
$p\varpi_K^r\cO_K=\varpi_K^{r+e}\cO_K$.

\begin{lemma}\label{app:lem-multiplicity-recurrence}
For $r\geq e+2$ and a nonzero finitely generated
$\Lambda$-module $N$,
\[
 \delta(\mathcal F(N))=\delta(N),\qquad
 \mu_{U_r}(\mathcal F(N))\geq q\,\mu_{U_{r-1}}(N).
\]
\end{lemma}

\begin{proof}
The group $U_r$ is normal in $I_1/Z_1$ and is contained in
$\mathcal K_1/Z_1$. Thus $\mathcal F(N)|_{U_r}$ is a direct
sum of $[I_1:\mathcal K_1]=q$ copies of $N$ with twisted action.
Precisely, for a representative $g\in I_1/\mathcal K_1$, the
action on the corresponding copy is
\[
 u\star_g v=(g^{-1}ug)\cdot_t v=(t^{-1}g^{-1}ugt)v,
 \qquad u\in U_r,\ v\in N.
\]
Since $U_r$ is normal in $I_1/Z_1$, conjugation by $g$ preserves
$U_r$ and its augmentation filtration. Since $t^{-1}U_rt$ is
open in $I_1/Z_1$, all the $q$ summands have canonical dimension
$\delta(N)$. Moreover, every summand has multiplicity
$\mu_{t^{-1}U_rt}(N)$, hence
\[
 \mu_{U_r}(\mathcal F(N))=q\,\mu_{t^{-1}U_rt}(N).
\]
Finally, one checks directly that $t^{-1}U_rt\subseteq U_{r-1}$. Both groups are uniform,
so Lemma~\ref{app:lem-multiplicity-scale} gives the inequality.
\end{proof}

\subsection{Proof of the main theorem}

Recall $M=(I^+/C_1)^\vee$, $Q=(C_0/C_1)^\vee$ and the short
exact sequence from \eqref{app:eq-exact}:
\[
 0\longrightarrow\mathcal F(M)\longrightarrow M
   \longrightarrow Q\longrightarrow0.
\]

\begin{proposition}\label{app:prop-dimension-bound}
If $M\ne0$, then $\delta(M)\geq n$. Moreover, if equality holds,
then $\delta(Q)<n$.
\end{proposition}

\begin{proof}
Put $t=\delta(M)$. For $r\geq e+1$, write
\[
 a_r=\mu_{U_r}(M)>0,\qquad b_r=\mu_{U_r,t}(Q)\geq0.
\]
For $r\geq e+2$, the above short exact sequence and
Lemma~\ref{app:lem-multiplicity-recurrence} imply
\begin{equation}\label{app:eq-multiplicity-recurrence}
 a_r=\mu_{U_r}(\mathcal F(M))+b_r\geq qa_{r-1}+b_r.
\end{equation}
Choose $m\geq 2e+1$ and apply this inequality to
$r\in\{m,m-1,\ldots,m-(e-1)\}$. We obtain
\[
 a_m\geq q^e a_{m-e}+\sum_{j=0}^{e-1}q^j b_{m-j}.
\]
On the other hand, \eqref{app:eq-congruence-powers} and
Lemma~\ref{app:lem-multiplicity-scale} give $a_m=p^t a_{m-e}$.
Since $q^e=p^{ef}=p^n$, it follows that
\begin{equation}\label{app:eq-key-inequality}
 (p^t-p^n)a_{m-e}\geq\sum_{j=0}^{e-1}q^j b_{m-j}\geq0.
\end{equation}
Thus $t\geq n$. If $t=n$, then $b_m=0$, which means
$\delta(Q)<n$.
\end{proof}

Finally we can prove the main result of this appendix.

\begin{proof}[Proof of Theorem~\ref{app:thm-canonical-diagram}]
Since $\pi$ is supersingular and $K\ne\Q_p$, it is not finitely
presented by the main results of \cite{Sch15} and \cite{Wu21}.
The finite-presentation criterion recalled above implies that
$C_1$ is infinite-dimensional. In particular, $\pi$ is
infinite-dimensional, and Corollary~\ref{cor:main-gk-dimension}
gives $\GKdim(\pi)=n$.

The inclusion $C_1\subseteq\pi$ induces a continuous surjection
$\pi^\vee\twoheadrightarrow C_1^\vee$ of $\F[[H]]$-modules,
since $H\subseteq I$. The latter module is finitely generated:
$C_1|_H$ is an admissible subrepresentation of $\pi|_H$.
Give it the quotient filtration. {This is a complete separated
good filtration, and its associated graded module is a quotient of
$\gr_\omega(\pi^\vee)$. } Hence
\[
 J^a\gr_\omega(C_1^\vee)=0.
\]
Theorem~\ref{thm:main-lambda} therefore applies to the
infinite-dimensional module $C_1^\vee$, giving
$\GKdim(C_1)=n$.

Next we prove $\GKdim(C_0/C_1)<n$. There are two cases.
\begin{itemize}
\item If $M=0$, then $I^+=C_1$, and $\Pi$-stability of
$C_1$ gives $I^-=C_1$ as well. Therefore
$\pi=I^++I^-=C_1=C_0$ and $C_0/C_1=0$.
\item If $M\ne0$, put $t=\delta(M)$. Then
$t=\GKdim(I^+/C_1)\leq\GKdim(\pi)=n$.
Combined with Proposition~\ref{app:prop-dimension-bound}, this
gives $t=n$ and $\GKdim(C_0/C_1)<n$.
\end{itemize}

Using Proposition~\ref{app:prop-purity}(iii), we deduce
\[
 W:=P_{n,I}(C_1)=P_{n,I}(C_0)\ne0.
\]
By part (ii), $W=P_{n,\mathcal K}(C_0)$, so $W$ is
$\mathcal K$-stable. On the other hand, $\Pi C_1=C_1$ and
$\Pi$ normalizes $I$, so uniqueness of $P_{n,I}(C_1)$ shows
that $W$ is $\Pi$-stable. Since $G$ is generated by
$\mathcal K$ and $\Pi$, irreducibility gives $W=\pi$.
As $W\subseteq C_1\subseteq C_0\subseteq\pi$, the result follows.
\end{proof}

\end{cbunit}


\end{document}